\documentclass[12pt]{article} 

\usepackage[utf8]{inputenc} 
\usepackage{amsmath,amsthm,verbatim,amssymb,amsfonts,amscd,graphicx}
\usepackage[colorlinks,citecolor=blue,urlcolor=blue]{hyperref}

\usepackage{geometry} 
\usepackage{graphicx} 

\usepackage{dsfont}

\usepackage{booktabs} 
\usepackage{array} 
\usepackage{verbatim} 
\usepackage{subfig} 
\usepackage{amsmath}
\usepackage[utf8]{inputenc}
\usepackage[english]{babel}
\usepackage{csquotes} 

\usepackage{algorithm2e} 

\usepackage{bm}

\usepackage[symbol]{footmisc}

\usepackage{fancyhdr} 
\usepackage{sectsty}
\allsectionsfont{\sffamily\mdseries\upshape} 

\usepackage[nottoc,notlof,notlot]{tocbibind} 
\usepackage{natbib}
\usepackage{bbm}
\usepackage[titles,subfigure]{tocloft} 
\usepackage{comment}
\usepackage{stmaryrd} 
\usepackage{enumitem} 

\newtheorem{theorem}{Theorem}[section]

\newtheorem{corollary}{Corollary}[section]

\newcommand{\qmq}[1]{\quad\mbox{#1}\quad}

\newcommand{\qm}[1]{\quad\mbox{#1}}
\newcommand{\abs}[1]{\left|{#1}\right|}
\newcommand{\wtilde}[1]{\widetilde{#1}}
\newcommand{\what}[1]{\widehat{#1}}
 
\newcommand*{\Var}[1]{\mathrm{Var} \left [ {#1} \right ]} 

\newcommand*{\E}[1]{\mathbb{E} \left [ {#1} \right ]}

\newcommand*{\R}{\mathbb{R}} 
\newcommand*{\diag}[1]{\mathrm{diag}({#1})}
\newcommand{\mA}{\mathcal{A}}

\newcommand{\mI}{\mathcal{I}}

\newcommand{\mN}{\mathcal{N}}

\newcommand{\To}{\rightarrow}
\newcommand{\eps}{\varepsilon}
\newcommand{\var}{\mbox{Var}}
\newcommand{\cov}{\mbox{Cov}}
\newcommand{\sd}{\mbox{SD}}

\title{Non-asymptotic bounds for quasi-MLE, misspecified models, and dependence under group sequential sampling}
\date{Last modified 18.May.26}
\author{\textsc{Julian Aronowitz}$^\dag$ and \textsc{Jay Bartroff}$^*$\\
\small{$^\dag$Google Inc., New York City, New York, USA}\\
\small{$^*$Department of Statistics and Data Sciences, University of Texas at Austin}\\
\small{Austin, Texas, USA}\\
}

\begin{document}
\maketitle

\abstract{We derive asymptotic multivariate normal limits and explicit non-asymptotic normal approximation bounds for group sequential quasi-maximum likelihood estimators under possible model misspecification and within-group dependence. The bounds, obtained using Stein's method, have known constants and apply to a class of dependent-data estimating problems in which the likelihood used for estimation may differ from the true data-generating mechanism. We compute the limiting covariance structure and finite-sample bound explicitly for a Poisson generalized linear mixed model with random group effects and illustrate the results using data from an epilepsy clinical trial.}

\section{Introduction and background}\label{ch:intro}

Group sequential analysis is a powerful statistical framework in which a fixed sample size is not pre-determined before conducting a study, but rather data is collected and analyzed in groups until the conditions of a pre-defined stopping criterion are met.  Group sequential methods are the dominant statistical methodology in many modern clinical trials \citep{sequential_book,Bartroff13} where they can provide ethical, administrative, and economic benefits, among other areas of application including
online A/B testing \citep{Kohavi20}, genomics and high-throughput biology \citep{Lin22}, and machine and reinforcement learning \citep{Powell22}. Critical to providing these benefits is the early stopping criterion, which in the most commonly used parametric settings is based on   maximum likelihood estimators (MLEs), repeatedly computed at each interim analysis, and thus whose joint distribution determines the operating characteristics of the trial.

Suppose sequential observations $Y_1, \dots , Y_n \in \R^t$ are divided into $K$ groups with $n_k$ denoting the number of observations up to and including group $k = 1, \dots, K$, and $n_K=n$. If $\what{\theta}_{n_k} =  \what{\theta}_{n_k}(Y_1, Y_2, \dots, Y_{n_k}) \in \R^{d}$ denotes the MLE of the parameter vector $\theta \in \R^d$ of interest based on the observations in the first $k$ groups, then the distribution of the concatenated vector
\begin{equation}\label{th.K.def}
\what{\theta}^K = [\what{\theta}_{n_1}; \what{\theta}_{n_2}; \dots; \what{\theta}_{n_K}] \in \mathbb{R}^{dK} 
\end{equation} 
is the key object in parametric group sequential methods. In \eqref{th.K.def} and throughout this paper we use the semicolon notation
\begin{equation}\label{transp.not}
[v_1;\ldots;v_K]\in\mathbb{R}^{dK}
\end{equation}
to denote the concatenation of column vectors $v_1,\ldots, v_K\in\mathbb{R}^d$ into a single  $\mathbb{R}^{dK}$-valued column vector.

Under the assumption of independent but not necessarily identically distributed observations,  \cite{jt} showed that $\what{\theta}^K$ is asymptotically multivariate normal under suitable regularity conditions, and \citet{Aronowitz25} provide non-asymptotic bounds to this normal limit. Although these results do not require identically distributed observations~$Y_i$, they do rely on independence of the $Y_i$ which has been a barrier to addressing statistical situations with dependence between observations such as correlated observations, random effects and mixed models, regularized estimators, and possibly misspecified models.  About this, \citet[][p.~1335-6]{jt} write, 

\begin{displayquote}
	\emph{``The general theory can also be extended to studies with correlated observations by replacing the distribution $f_i(Y_i; \theta)$ with the conditional distribution of $Y_i$ given $Y_1, \dots, Y_{i-1}$ and defining efficient scores\ldots The conclusion
    \ldots remains valid if it is permissible to apply a WLLN to sums of conditional information and a central limit theorem to sums of conditional scores to deduce the conditions required.''}
\end{displayquote}
We are not aware of such an extension to correlated observations having been  carried out in the literature  in any general setting.
Correlated observations are considered by \citet[][p.~67]{sequential_book} but only in the context of normal linear models; in this case the group sequential MLE is exactly multivariate normal and no asymptotic theory is required. A survey of results regarding group sequential analysis of specific longitudinal models is presented in \citet[p.~232-233]{sequential_book}.

The purpose of this paper is to develop a theory covering asymptotic normality results as well as nonasymptotic bounds to these limits for group sequential MLEs in a setting that allows dependence within groups. Examples of this include cluster-randomized trials and multicenter studies with site effects; in these settings, the i.i.d.\ assumption is inappropriate at the observation level, but independence may still be plausible at the group (e.g., cluster or site) level.  In doing so we additionally address misspecified models and so-called quasi maximum likelihood estimators (QMLEs). In the statistics literature, the term ``quasi-likelihood'' has been used in related yet distinct meanings. In the setting of \citet{mccullagh1983} and \citet{mcculloch2011_book}, a quasi-likelihood is a function that describes the relationship between the mean and variance of observed random variables without specifying the full distribution of the data generating process. Optimizing this ``variance'' function produces an estimator that has similar properties to the MLE. The second usage of quasi-likelihood is the context originated by \citet{white1982,white1996_book}, who characterized a quasi-likelihood function as one essentially equivalent to the true likelihood except for the fact that there may not be a parameter in the parameter space for which the quasi-likelihood correctly specifies the data generating process of interest, a situation known as model misspecification. The QMLE is then the parameter in the parameter space that maximizes the quasi-likelihood given the observed data. This paper is concerned with this latter meaning of quasi-likelihood.

\citet{white1996_book} showed asymptotic normality of QMLEs in the fixed sample setting: Under regularity conditions the QMLE~$\what{\theta}$ is consistent for some $\theta^*$ that minimizes the Kullback-Leibler (KL) divergence (see \eqref{KL.def}) from the true data generating process~$g(y)$ to the quasi-likelihood~$f(y|\theta)$, and that the QMLE is asymptotically normal around $\theta^*$.  Our Theorem~\ref{thm:qmle.seq.asymp.norm} extends and generalizes this result to find the asymptotic distribution of QMLEs in the group sequential setting, and Theorem~\ref{thm:qmle_non_asymptotic} provides  a finite-sample bound, with explicit constants, for QMLEs based on dependent data. By taking the specified density~$f$ to assume independence while the true density~$g$ does not, we construct misspecified although consistent QMLE estimators that are asymptotically normal, and Theorem~\ref{thm:qmle_non_asymptotic} quantifies the effect of this misspecified dependence on the normal limit.  To give an example of how one might apply Theorem~\ref{thm:qmle_non_asymptotic}, in Section~\ref{sec:PoissGLM} we apply these results to the group sequential analysis of a Poisson generalized linear mixed model, and in Section~\ref{sec:data} compute the bound numerically using data from an epilepsy clinical trial.

The closest antecedent in the non-asymptotic MLE literature is \citet{AR_bound}, who obtained explicit Stein bounds for fixed-sample MLEs. Our results extend this line in three directions that arise naturally in sequential statistical problems: vector-valued group sequential estimators, model misspecification through QMLEs, and dependence within sampling groups. This places the paper at the intersection of Stein-type normal approximation, sequential inference, and asymptotic theory for dependent statistical models.

\subsection{Other related literature}
The non-asymptotic bound for the univariate MLE first developed by \citet{AR_bound} was expanded to the multivariate case by \citet{mvn}. In \citet{delta} and \citet{Anastasiou20} the bounds for the univariate and multivariate cases, respectively, were sharpened and simplified under the additional assumption that the MLE follows a certain additive form.  \citet{mle_bound_m_dependence} considers data that is $m$-dependent for identically distributed scalar parameter MLEs. However, our consideration of QMLE and misspecified models is new in this literature, and our results extend the recent work of  \citet{Aronowitz25} for correctly specified  models.  All of these results mentioned so far utilize Stein's \citeyearpar{Stein72} method  for the normal distribution, for which we refer the reader to \citet{larrybook} for an introduction.  Results not based on Stein's method include Pinelis' \citeyearpar{kol_mle_bound} optimal order bound for univariate MLEs in terms of the Kolmogorov distance, and bounds derived using characteristic functions \citep[see][]{ulyanov1979, ulyanov1986, ulyanov1987}. These results require independence and apply only to the univariate setting, and thus are not considered for extension to the setting considered here. The convergence of statistics after stopping a group sequential trials has been considered by B.~Berckmoes and coauthors \citep{Berckmoes18,Berckmoes18b,Berckmoes20}. Regarding dependence, our results require a slightly stronger assumption than conjectured in the quote above, namely we allow dependence within groups but require independence between groups.  Under this condition, the ``independent increments'' structure of the group sequential analysis is preserved \citep[see][]{Spiessens00,Kim20,Parast24}.

\section{QMLEs and dependent data}\label{ch:qmle_and_dependent_data}

In Section~\ref{sec:qmle.intro} we define the quasi-maximum likelihood estimator and give White's~\citeyearpar{white1996_book} result for its limiting normality.  In Section~\ref{sec:results_for_misspecified} we extend this to the group sequential setting, and then derive non-asymptotic bounds for the group sequential limit in Section~\ref{sec:GS.bds}.

\subsection{Asymptotic distribution: Notation and White's~\citeyearpar{white1996_book} fixed-sample result}\label{sec:qmle.intro}

We begin in the fixed-sample setting in which the sample size denoted by $n$, which will approach infinity in the asymptotic results below.   Letting $Y^n=(Y_1,\ldots, Y_n)$, $Y^0=\emptyset$,  and similar for lowercase arguments~$y^n$, we let  $g^n(y^n)$ denote the joint density of the true, unknown data generating process, and let 
\begin{equation}\label{fn.def}
f^n(y^n |\theta) =\prod_{i = 1}^{n} f_i(y_i|y^{i-1}, \theta),\quad  \theta \in \Theta \subseteq \R^d,
\end{equation}
denote a postulated family of parametric models for $Y^n$. Mutual independence of $Y_1,\ldots, Y_n$ is \emph{not} assumed in $g^n$ or $f^n$, as the conditional notation in \eqref{fn.def} reflects. We treat a particular $f^n(\cdot|\theta)$ as an approximation of the true $g^n$ whose accuracy is measured by the KL divergence
\begin{equation}\label{KL.def}
\mathbb{E} \left[\log\left(\frac{g^n(Y^n)}{f^n(Y^n|\theta)}\right)\right]
\end{equation}
between $g^n$ and $f^n(\cdot|\theta)$. In \eqref{KL.def} and  throughout, the expectation is taken with respect to $Y^n\sim g^n$, the true data generating process. Minimizing the KL divergence~\eqref{KL.def} with respect to $\theta$ is equivalent to choosing the $\theta$ that maximizes the expected log-likelihood, i.e.,
\begin{equation} \label{eq:beta^*}
\theta_n^* = \arg \max_{\theta \in \Theta} \E{\log f^n(Y^n|\theta)}.
\end{equation}
If there exists $\theta_0 \in \Theta$ such that $f^n(Y^n|\theta_0) = g^n(Y^n)$ a.s.\ then we say that $f^n$ is \textit{correctly specified}. In this case, the minimizer~$\theta_n^*$ of the KL divergence  is equal to the true parameter value~$\theta_0$. Since the expectation~\eqref{eq:beta^*}
is with respect to the unknown $g^n$, computing $\theta_n^*$ directly is not an option available to the statistician. A natural alternative is the  empirical version of \eqref{eq:beta^*},
\begin{equation}\label{QMLE.def}
\what{\theta}_n = \arg \max_{\theta \in \Theta} n^{-1} \log f^n(Y^n|\theta) = \arg \max_{\theta \in \Theta} n^{-1}\sum_{i = 1}^n\log f_i(Y_i |Y^{i-1}, \theta),
\end{equation}
which we call the \textit{quasi-maximum likelihood estimator (QMLE)} of $\theta$.

\citet[p.~92]{white1996_book} showed that, under proper regularity conditions, the QMLE exhibits similar asymptotic normality properties to the MLE, but centered at $\theta_n^*$. Specifically, that 
\begin{equation}\label{white.CLT}
\sqrt{n} J^{-1/2} (\what{\theta}_n - \theta_n^*) \overset{\mathcal{D}}{\longrightarrow} \mN(0,\mI_d)\qmq{as}n\To\infty,
\end{equation} denoting convergence in distribution to the $d$-dimensional multivariate normal with mean vector~$0$ and covariance matrix the $d\times d$ identity  matrix~$\mI_d$, for a certain matrix $J=J(\theta_n^*)$ whose group sequential generalization is given in Theorem~\ref{thm:qmle.seq.asymp.norm}. The needed regularity conditions for \eqref{white.CLT} are listed in Appendix~\ref{sec:reg.conds}, and are the same as those we will require for our group sequential generalization of \eqref{white.CLT} in Theorem~\ref{thm:qmle.seq.asymp.norm}.

\subsection{Asymptotic distribution of group sequential QMLE for misspecified models}\label{sec:results_for_misspecified}

In Theorem~\ref{thm:qmle.seq.asymp.norm} we find the asymptotic joint distribution of the group sequential QMLE~$\what{\theta}^K$, given by \eqref{th.K.def},  under a possibly misspecified model, generalizing \eqref{white.CLT} to the group sequential setting. The main distributional assumptions are the same ones needed for \eqref{white.CLT} but, to simplify the theorem's  presentation, we make the additional assumption of 
mutual independence between, but not necessarily within, groups; a more general result, allowing some dependence between groups as well, will hold as long as the LLN, consistency, and  marginal asymptotic normality used hold for the conditional score functions. 

For the remainder of the paper we let $[j]=\{1,2,\ldots, j\}$, and we add the following notation to the group sequential notation introduced in Section~\ref{ch:intro}.  Let \begin{equation}\label{theta.star.GS}
\theta^{*K} = [\theta_{n_K}^{*}; \theta_{n_K}^{*}; \dots; \theta_{n_K}^{*}] \in \R^{q} \qm{where $q=dK$,}
\end{equation} 
recalling our notation~\eqref{transp.not}. With 
$n_0=0$, let 
\begin{equation}\label{Gk.def}
G_k=\{n_{k-1}+1, n_{k-1}+2,\ldots, n_k\},\quad k\in[K],  
\end{equation} denote the indices of the $k$th group, let $Y^{(k)}=\{Y_i:\; i\in G_k\}$, $k\in[K]$, denote the observations in the $k$th group only,  
and let
\begin{equation}\label{score.def}
S(i,\theta) = \nabla \log f_i(Y_i|Y^{i-1}, \theta) \in \R^{d},\quad i\in[n_K],
\end{equation} 
be the conditional score function, a column vector. 
The asymptotic regime we consider in Theorem~\ref{thm:qmle.seq.asymp.norm} is that in which, for every group $k\in[K]$, the cumulative group size diverges~$n_k\To\infty$ in such a way that the fraction $(n_k-n_{k-1})/n_K$ of observations in the $k$th group approaches some constant in $(0,1)$, possibly depending on $k$.  We denote this asymptotic regime by ``$\{n_k\}\To\infty$'' and, to ease notation, refrain from indexing some objects by $n_1,\ldots, n_K$ which depend on them, such the matrices in \eqref{eq:Jn_misspecified_definition}-\eqref{Ik.def}. We use $\wedge$ and $\vee$ for min and max, respectively. Finally, we will utilize matrices and vectors of size $q=dK$ with block structure of size $d$, denoted as follows:  For $v\in\mathbb{R}^q$ let $v_{[i]}\in\mathbb{R}^d$ denote the $i$th sub-vector (or ``$d$-block''), for $i\in[K]$, and for a $q\times q$ matrix~$B$ let $B_{[i][j]} \in \R^{d \times d}$ denote the $(i,j)$ block, for $i,j\in[K]$. 

\begin{theorem}\label{thm:qmle.seq.asymp.norm}
Assume that regularity conditions~\ref{reg:1st}-\ref{reg:CLT} in Appendix~\ref{sec:reg.conds} hold, and that $Y^{(1)}, Y^{(2)}, \ldots, Y^{(K)}$ are mutually independent. Then
\begin{equation*}
\sqrt{n_K}J^{-1/2}(\what{\theta}^K - \theta^{*K}) \overset{\mathcal{D}}{\longrightarrow} \mN(0,\mI_q)\qmq{as}\{n_k\}\To\infty
\end{equation*}
where, for $j,k\in[K]$,
\begin{align}
J_{[j][k]} &= H_j(\theta_{n_K}^*)^{-1}I_{j\wedge k}(\theta_{n_K}^*)H_k(\theta_{n_K}^*)^{-1},\label{eq:Jn_misspecified_definition}\\
H_k(\theta) &= -\E{ \nabla \frac{1}{n_K} \sum_{i = 1}^{n_k}S(i,\theta)}, \label{Hk.def}\\
I_k(\theta) &= \Var{\frac{1}{\sqrt{n_K}}\sum_{i = 1}^{n_k} S(i,\theta)},\label{Ik.def}
\end{align} and the score function~$S(i,\theta)$ is given by \eqref{score.def}.
\end{theorem}

The theorem's proof is in Appendix~\ref{sec:GS.asymp.pf}. 

In the generality of the theorem which allows model misspecification, the Hessian and information matrices \eqref{Hk.def}-\eqref{Ik.def} may differ. When the model is correctly specified in the sense of the discussion following \eqref{eq:beta^*}, then \begin{equation}\label{IME}
H_k(\theta_0)=I_k(\theta_0),
\end{equation}
known as the \textit{Information Matrix Equality} \citep{Greene08,Hayashi11}, and  Theorem~\ref{thm:qmle.seq.asymp.norm} reduces to the correctly-specified special case in \citet{Aronowitz25} for group sequential sampling, while the fixed-sample analog -- White's result \eqref{white.CLT} -- reduces to the classical limiting normality of the MLE \citep[see][]{Cox00}.

A natural question in applying the multivariate normal limit in Theorem~\ref{thm:qmle.seq.asymp.norm} is how quickly it is reached. Next we address this by finding non-asymptotic bounds to these limits.

\subsection{Non-asymptotic bounds}\label{sec:GS.bds}

In Theorem~\ref{thm:qmle_non_asymptotic} we give an upper bound on $\abs{\mathbb{E} [h(X)] - \mathbb{E} [ h(Z)]}$, where $X$ is the normalized group sequential QMLE, $Z$ is a $q$-dimensional standard normal,  and $h$ is a test function in the smooth function class~$\mathcal{H}$, given in \eqref{H}, as follows. For any $3$ times differentiable function $h: \R^q \rightarrow \R$ we denote $|h|_0 := \sup \abs{h}$ and 
$$|h|_i = \sup_{\sum_j \alpha_j=i} \left|\frac{\partial^{i}h}{\partial x_1^{\alpha_1}\cdots \partial x_q^{\alpha_q}}\right|,\quad i=1,2,3,$$ where the supremum is taken over all multi-indices $(\alpha_1,\ldots, \alpha_q)$ summing to $i$, as well as arguments to the function.  The function class considered is then
\begin{equation}\label{H}
\mathcal{H} = \{h: \R^q \rightarrow \R\;\mbox{$3$ times differentiable with $|h|_i<\infty$, $i=0,1,2,3$}\}.
\end{equation}

In order to state the theorem we must introduce some additional notation. Extend the notation \eqref{score.def} by letting 
\begin{equation}\label{S.set.notation}
S(s, \theta) = \sum_{i\in s} S(i,\theta) \qmq{for any set $s\subseteq [n_K]$.}
\end{equation}
Let
\begin{equation}\label{Itilde.def}
\wtilde{I}_k(\theta)= \Var{ n_K^{-1/2} S(G_k,\theta)},\quad k\in[K],
\end{equation} the corresponding group-wise information matrices,
\begin{equation}\label{non.asymp.c.def}
\tau= \max_{k \in [K],\; i,j\in[d]} \abs{\wtilde{I}_k(\theta_{n_K}^*)_{ij}^{-1/2}},
\end{equation}
and
\begin{equation}\label{xi.def}
\xi_i=n_K^{-1/2} S(i,\theta_{n_K}^*), \quad i\in[n_K],
\end{equation} with $\xi_i'$ an independent copy of $\xi_i$. Let 
\begin{equation}\label{Q.def}
Q_{ki} = (\what{\theta}_{n_k} - \theta_{n_K}^*)_i,\quad k\in [K], i\in [d],\qmq{and} Q=\max_{k\in [K],\; i\in [d]} \left|Q_{ki}\right|.
\end{equation} 
Finally, let $\sd(\cdot)$ denote $\sqrt{\var(\cdot)}$.

\begin{theorem}\label{thm:qmle_non_asymptotic}
Assume that the conditions of Theorem~\ref{thm:qmle.seq.asymp.norm} are satisfied, and additionally that for any value of $\theta_{n_K}^*$ there exists $0<\eps = \eps(\theta_{n_K}^*)$ and functions $M^k_{iuj}(y)$, $i, u, j \in [d]$, $k\in [K]$, such that
\begin{equation}\label{M.fcn.bd}
\abs{ \frac{\partial^3 \log f^{n_k}(y^{n_k}|\theta)}{\partial \theta_i \partial \theta_u \partial \theta_j}  } \leq M^k_{iuj}(y^{n_k})
\end{equation}
for all $\theta \in \Theta$ with $|\theta_j - \theta_{n_K j}^*| < \eps$ for all $j \in [d]$, and
\begin{equation}\label{EM.fcn.bd}
\mathbb{E} \left[ \left. \left(  M^k_{iuj}(Y^{n_k}) \right)^2 \right| Q < \eps \right] < \infty.
\end{equation}
Then for $Z \sim \mN(0, \mI_q)$, $h \in \mathcal{H}$ given by \eqref{H}, and $J$ defined by \eqref{eq:Jn_misspecified_definition}, we have 
\begin{multline}\label{eq:qmle_non_asymptotic_statement}
\abs{\mathbb{E} h(\sqrt{n_K} J^{-1/2} (\what{\theta}^K - \theta^{*K})) - \mathbb{E} h(Z)} \\ 
\le \frac{|h|_1}{\sqrt{n_K}} R_1 + \frac{ q^2 \tau^2 |h|_2}{4} R_2 + \frac{  q^3 \tau^3 |h|_3}{12} R_3 + 2|h|_0\frac{\mathbb{E}Q}{\eps},
\end{multline}
where
\begin{align}
R_1 &= \sum_{k_1=1}^K \sum_{k_2=1\vee (k_1-1)}^{k_1} \sum_{i,j=1}^d  \left| \wtilde{I}_{k_1}(\theta_{n_K}^*)_{ij}^{-1/2}\right|  \label{Rem1.def}\\
&\times \left\{ \sum_{l = 1}^d \sd\left[\nabla S([n_{k_2}],\theta_{n_K}^*)_{jl}\right] \sqrt{\mathbb{E}[Q_{k_2 l}^2]} + \left| \mathbb{E}S([n_{k_2}], \theta_{n_K}^*)_j \right|\right.\label{R1.def.L2}\\
&+ \left.\frac{1}{2} \sum_{l,l' = 1}^d \sqrt{\mathbb{E}[(Q_{k_2l}Q_{k_2l'})^2\wedge \eps^4] \; \mathbb{E} [ (  M_{jll'}^{k_2}(Y^{n_{k_2}}) )^2 | Q< \eps ]}\right\},\label{R1.def.L3}\\
R_2 &= \sum_{k = 1}^{K} \left\{ \sum_{j = 1}^d \sd\sum_{v\in G_k}\xi_{vj}(\xi_{vj}-2\mathbb{E} \xi_{vj}) \right. \nonumber\\
&+\left.  \sum_{i<j} \sd \sum_{v \in G_k}\left(\xi_{vi}\xi_{vj}-\xi_{vj}\mathbb{E}\xi_{vi} - \xi_{vi}\mathbb{E}\xi_{vj}\right) \right\},\qm{and}\nonumber\\
R_3 &= \sum_{i = 1}^{n_K} \mathbb{E} \left( \sum_{j = 1}^d \abs{\xi'_{ij} - \xi_{ij}}  \right)^3.\nonumber
	\end{align}
\end{theorem}

The proof of Theorem~\ref{thm:qmle_non_asymptotic} is in Appendix~\ref{sec:qmle_non_asymptotic}. The bound~\eqref{eq:qmle_non_asymptotic_statement} has components similar to  the non-asymptotic bound for correctly specified models in \citet[][Theorem~3.1]{Aronowitz25} but is more complex due to the possible dependence within groups and the fact that, in Theorem~\ref{thm:qmle_non_asymptotic}, the estimator~$\what{\theta}_{n_K}$ concentrates around  $\theta_{n_K}^*$ rather than $\theta_0$ and, unlike the latter, $\mathbb{E}S([n_K],\theta_{n_K}^*)$ does not necessarily vanish.

The key tool to finding the bound~\eqref{eq:qmle_non_asymptotic_statement} is a result of \citet{RR_exchangeable_pair_bound}, stated  in Appendix~\ref{app:Stein.bg}, which involves  the first $3$ derivatives of the test function~$h$.  \citet[][Proposition 2.1]{Gaunt16} found new bounds that require one fewer
derivative of $h$, and  \citet[][Theorem~3.6]{Gaunt23} have an improved result which permits a
version of this bound that depends on $|h|_1$ and $|h|_2$ but not on $|h|_0$. It may be possible to produce a version of our bound~\eqref{eq:qmle_non_asymptotic_statement} with these relaxations, at the cost of an increase of the bound's order of the dimension~$q$.

The next section is an in-depth application of Theorems~\ref{thm:qmle.seq.asymp.norm} and \ref{thm:qmle_non_asymptotic} to mixed Poisson regression models that   illustrates the within-group dependence, quasi-likelihood, and model misspecification allowed by those theorems.

\section{Poisson regression with random group effects}\label{sec:PoissGLM}

We apply the above theory to the case of Poisson regression with random group effects. After setting up the model's notation in Section~\ref{sec:pois.setup},  the parameters of the group sequential QMLE's asymptotic distribution are found in Section~\ref{sec:pois.asymp}, and the corresponding non-asymptotic bound is calculated in Section~\ref{sec:pois.non-asymp}.   These quantities are calculated for an epilepsy clinical trial data set in Section~\ref{sec:data}.

\subsection{Set up}\label{sec:pois.setup}
Consider the regression setting where the $i$th observation~$Y_i$ has associated covariate vector~$x_i\in\mathbb{R}^d$, and the $k$th group has associated random  effect~$U_k$. Adopting the Poisson generalized linear mixed model with the canonical link function, and letting $p(\cdot|\lambda)$ denote the Poisson density with mean~$\lambda$, the observations~$Y_i$ in the $k$th group are conditionally independent with
\begin{equation}\label{eq:poisson_glmm}
Y_i \mid U_k \sim p(\cdot|\exp(x_i^T\theta +U_k)) \qmq{for} i \in G_k.
\end{equation}
We further assume that the random effects~$U,U_1,U_2,\ldots, U_K$ are i.i.d.\ with density $\pi(\cdot|\psi)$, known up to parameter~$\psi$, with respect to some measure~$\nu$, and that the covariate vectors $x_1,\ldots, x_{n_K}$ have full rank~$d$. Thus the unknown parameters of the model are $(\theta, \psi)$, although in what follows we focus on estimation of $\theta$ and treat $\psi$ as a nuisance parameter, which we sometimes suppress in notation.

We consider a QMLE set up in which the postulated density~$f$ ignores the within-group dependence introduced by the random group effects in \eqref{eq:poisson_glmm}.  In practice, use of this $f$ may be a choice to simplify the analysis, or may be arrived at out of ignorance of the group effect.  In the following we let $\prod_{k,i}$ be shorthand for $\prod_{k\in[K]}\prod_{i\in G_k}$, and similar for $\sum_{k,i}$.  Straightforward calculations give that the true density is given by
\begin{multline}\label{pois.g.def}
g^{n_K}(y^{n_K}|\theta) = \int\cdots\int \prod_{k\in [K]} \left(\prod_{i\in G_k} p(y_i|e^{x_i^T\theta+U_k})\right) \pi(u_k|\psi) \nu(du_k)\\
= \prod_{k\in [K]} \int  \left(\prod_{i\in G_k} \left[\left. \exp\left(-e^{x_i^T\theta +u}\right) e^{y_i(x_i^T\theta +u)}\right/ y_i!\right]\right) \pi(u|\psi) \nu(du)\\
= \prod_{k\in [K]} \int  \left(\prod_{i\in G_k} p(y_i|e^{x_i^T\theta}) \exp\left[e^{x_i^T\theta}(1-e^u)+y_iu\right]\right) \pi(u|\psi) \nu(du)\\
= \left(\prod_{k', i'} p(y_{i'}|e^{x_{i'}^T\theta})\right)\prod_{k\in[K]} \int  \exp\left[\sum_{i\in G_k} (e^{x_i^T\theta}(1-e^u)+y_iu)\right] \pi(u|\psi) \nu(du)\\
=f^{n_K}(y^{n_K}|\theta)\prod_{k\in[K]}\wtilde{g}_k
\end{multline}
where
\begin{align}
f^{n_K}(y^{n_K}|\theta) &= \prod_{k,i}    p(y_i|e^{x_i^T\theta}),\label{pois.f.def}\\
\wtilde{g}_k &=\exp\left(\sum_{i'\in G_k}e^{x_{i'}^T\theta}\right)\int \exp\left[\sum_{i\in G_k} (y_iu -e^{x_i^T\theta} e^u)\right] \pi(u|\psi) \nu(du).\label{pois.gk.part}
\end{align}
\citet[p.~226]{mcculloch2011_book} note that \eqref{pois.gk.part}, and hence \eqref{pois.g.def}, cannot be simplified further or evaluated in closed form, even if $U_k$ is normally distributed.  For just computing the MLE of $\theta$, numerical integration of \eqref{pois.g.def} is typically used; other techniques for estimating $\theta$ are discussed in \citet{glmm}.  But for more refined analysis such as that required to compute non-asymptotic bounds, the  within-group dependence due to the random effect makes direct analysis of \eqref{pois.g.def} difficult.  For example, the  resulting score function is not simply the sum of the scores of the observations.  For these reasons, we next proceed to consider the QMLE approach using $f(\cdot|\theta)$ in \eqref{pois.f.def} as the postulated density family, assuming the true density is $g(\cdot|\theta_0)$ given by \eqref{pois.g.def} for some true,  unknown value $\theta_0$.

\subsection{Asymptotic distribution of QMLE for Poisson GLMM}\label{sec:pois.asymp}
In this section we apply  Theorem~\ref{thm:qmle.seq.asymp.norm} to the Poisson GLMM to obtain the group sequential asymptotic distribution.  Our results apply to the fixed sample-$n$ set up as well by taking $K=1$ and $n_K=n$, which has not appeared in the literature as far as we know.  While many references \citep[e.g.,][]	{Breslow93,McCulloch12,Wang22} discuss the Poisson GLMM, none of these compute explicitly the Hessian and information matrices in our setting.  In doing so we will see that the value~$\theta_{n_K}^*$ around which the asymptotic distribution concentrates  is the true value~$\theta_0$, scaled by the factor~$\mathbb{E} e^{U}$ determined by the distribution of the random effects.  In other words, if 
\begin{equation}\label{EU1}
\mathbb{E} e^{U}=1,   
\end{equation} then $\theta^*=\theta_0$. Thus, to simplify the notation in what follows, we assume \eqref{EU1} throughout. Otherwise, our analysis can be carried out with additional notational load incorporating the scaling factor of $\mathbb{E}e^{U}$ or, equivalently, fixing the offset in the Poisson regression model~\eqref{eq:poisson_glmm}.  After verifying $\theta_{n_K}^*=\theta_0$, we calculate the Hessian and information matrices~$H_k$ and $I_k$ in \eqref{Hk.def}-\eqref{Ik.def}, before applying the non-asymptotic bounds of Theorem~\ref{thm:qmle_non_asymptotic} in the next section.

Regarding the moments of $Y_i$, it is well-known \citep{Riordan37} that the $m$th moment ($m\ge 1$) of a Poisson random variable with mean~$\lambda$ is given by the polynomial in $\lambda$
\begin{equation}\label{pois.mom.gen}
\sum_{j=1}^m {m\brace j} \lambda^j,\qmq{where} {m\brace j}=\frac{1}{j!}\sum_{\ell=1}^j {j\choose \ell} (-1)^{j-\ell} \ell^m
\end{equation}
are the Stirling numbers of the second kind \citep[see][p.~125]{vanLint01}; special cases include ${m\brace 1}={m\brace m}=1$, and ${m\brace j}=0$ for $j>m$.  Applying this to $Y_i$, conditionally on $U_k$ for $i\in G_k$, we have
\begin{multline}
\mathbb{E}(Y_i^m) = \mathbb{E} \sum_{j=1}^m {m\brace j} (e^{x_i^T\theta_0+U_k})^j =\sum_{j=1}^m {m\brace j} e^{jx_i^T\theta_0}\mathbb{E}e^{jU_k}\\
=\sum_{j=1}^m {m\brace j} e^{jx_i^T\theta_0}M_U(j),  \label{pois.Y.mom}
\end{multline}
where $M_U(j) = \mathbb{E}e^{jU}$ is the moment generating function (MGF) of $U$. Using \eqref{pois.Y.mom} and $M_U(1)=1$,
\begin{align}
E(Y_i) &= e^{x_i^T\theta_0},\label{pois.EY}\\
E(Y_i^2)&= e^{x_i^T\theta_0} + e^{2x_i^T\theta_0}M_U(2).\label{pois.EY2}
\end{align}

Using \eqref{pois.EY}, recall from \eqref{eq:beta^*} that $\theta_{n_K}^*$ maximizes 
\begin{multline}\label{pois.Ef}
\mathbb{E}[\log f^{n_K}(Y^{n_K}; \theta)] = \sum_{k,i} \mathbb{E}[\log p(Y_i|e^{x_i^T\theta})] =\sum_{k,i} \left(-e^{x_i^T\theta}+x_i^T\theta \mathbb{E}Y_i\right)+C\\
=\sum_{k,i} \left(-e^{x_i^T\theta}+x_i^T\theta e^{x_i^T\theta_0}\right)+C,
\end{multline}
where $C$ does not depend on $\theta$. The gradient of \eqref{pois.Ef} with respect to $\theta$ is
$\sum_{k,i} (e^{x_i^T\theta_0}- e^{x_i^T\theta})x_i$, and it follows from the full-rank assumption on the $x_i$ that this vanishes if and only if $\theta=\theta_0$. Therefore, in what follows we use $\theta_0$ in place of $\theta_{n_K}^*$, and let 
\begin{equation}\label{pois.ei.def}
e_i=e^{x_i^T\theta_0},\quad i\in [n_K],  
\end{equation}
throughout this and the next section.

We have 
\begin{equation}\label{pois.S}
S(i,\theta) = \nabla \log p(Y_i|e^{x_i^T\theta})=(Y_i-e^{x_i^T\theta})x_i\qmq{and} \nabla S(i,\theta) = -e^{x_i^T\theta} x_ix_i^T,
\end{equation}
so the Hessian~\eqref{Hk.def} is
\begin{equation}\label{pois.H.t0}
H_k(\theta_0) =\frac{1}{n_K}\sum_{i=1}^{n_k} e_i x_ix_i^T. 
\end{equation}
To compute the cumulative information matrix~\eqref{Ik.def}, we first compute the group-wise version~\eqref{Itilde.def}.  Note that 
\begin{equation}\label{pois.ES=0}
\mathbb{E}[S(i,\theta_0)]=0
\end{equation}
which follows from \eqref{pois.EY} and \eqref{pois.S}. Then
\begin{multline}
n_K \wtilde{I}_k(\theta_0) =  \var\left[\sum_{i\in G_k} S(i,\theta_0)\right]=  \mathbb{E}\left[\sum_{i\in G_k} S(i,\theta_0)\right] \left[\sum_{j\in G_k} S(j,\theta_0)\right]^T \\
=  \sum_{i\in G_k} \mathbb{E}\left[S(i,\theta_0) S(i,\theta_0)^T\right] +\sum_{i\ne j\in G_k} \mathbb{E}\left[S(i,\theta_0)S(j,\theta_0)^T\right]\\
=  \sum_{i\in G_k} \mathbb{E}\left[(Y_i-e_i)^2\right]x_ix_i^T + \sum_{i\ne j\in G_k} \mathbb{E}\left[(Y_i-e_i)(Y_j-e_j)\right]x_ix_j^T\\
=  \sum_{i\in G_k} \var(Y_i) x_ix_i^T + \sum_{i\ne j\in G_k} \cov(Y_i,Y_j) x_ix_j^T.\label{pois.Itild.part}
\end{multline}
For $i\ne j\in G_k$, using conditional independence,
\begin{multline*}
\mathbb{E}[Y_iY_j] = \mathbb{E}[\mathbb{E}[Y_iY_j|U_k]] = \mathbb{E}[\mathbb{E}[Y_i|U_k] \mathbb{E}[Y_j|U_k]] = \mathbb{E}[(e_i e^{U_k})(e_j e^{U_k})]\\
= e_i e_j \mathbb{E}[e^{2U_k}] = e_i e_j M_U(2).
\end{multline*}
Using \eqref{pois.EY}, \eqref{pois.EY2}, and this last, 
\begin{align*}
\var(Y_i)&= e_i[1+e_i(M_U(2)-1)],\\
\cov(Y_i,Y_j)&= e_i e_j (M_U(2)-1)\qmq{for} i\ne j\in G_k.
\end{align*}
Plugging these into \eqref{pois.Itild.part} and writing $M_U(2)-1=\var (e^U)$, 
\begin{align}
\wtilde{I}_k(\theta_0) &=\frac{1}{n_K}\left\{\sum_{i\in G_k} e_i[1+e_i\var (e^U)] x_ix_i^T + \var (e^U) \sum_{i\ne j\in G_k} e_i e_j x_ix_j^T\right\},\label{pois.Itild.t0}\\
I_k(\theta_0)&=\sum_{k'\le k} \wtilde{I}_{k'}(\theta_0)\label{pois.I.t0}
\end{align} by independence between groups.

Comparing \eqref{pois.H.t0} and \eqref{pois.Itild.t0} we see that these correspond -- i.e., the Information Matrix Equality~\eqref{IME} holds -- when $\var (e^U)=0$, which with \eqref{EU1} imply that $U=0$ a.s. In this case the model~\eqref{eq:poisson_glmm} reduces to the postulated density~$f$ and is thus correctly specified.

\subsection{Non-asymptotic bounds}\label{sec:pois.non-asymp}
In this section we calculate the quantities  needed to apply the bounds in Theorem~\ref{thm:qmle_non_asymptotic}  to the Poisson GLMM \eqref{eq:poisson_glmm}. By the discussion following \eqref{pois.Ef}, we use $\theta_0$ in place of $\theta_{n_K}^*$ throughout.

Here, by \eqref{pois.S}, the variate~$\xi_i$ appearing in $R_2$ and $R_3$ is
$$\xi_i = n_K^{-1/2}S(i,\theta_0) = n_K^{-1/2} (Y_i-e_i)x_i$$ and $\mathbb{E}\xi_i=0$ by \eqref{pois.ES=0}.

We start with the second variance term in $R_2$, which we compute using the central moments of the Poisson conditionally:
\begin{multline*}
n_K^2 \mathbb{E}\left(\sum_{v\in G_k} \xi_{vi}\xi_{vj} \right)^2  \\
= \sum_{v\in G_k} x_{vi}^2 x_{vj}^2 \mathbb{E}(Y_v-e_v)^4+2\sum_{v<v'\in G_k} x_{vi}x_{vj}x_{v'i}x_{v'j} \mathbb{E}[(Y_v-e_v)^2 (Y_{v'}-e_{v'})^2]\\
= \sum_{v\in G_k} x_{vi}^2 x_{vj}^2 \mathbb{E}[3(e_ve^{U_k})^2+e_ve^{U_k}] +2\sum_{v<v'\in G_k} x_{vi}x_{vj}x_{v'i}x_{v'j} \mathbb{E}[(e_ve^{U_k})(e_{v'}e^{U_k})]\\
= \sum_{v\in G_k} x_{vi}^2 x_{vj}^2 [3 e_v^2M_U(2)+e_v] +2\sum_{v<v'\in G_k} x_{vi}x_{vj}x_{v'i}x_{v'j} e_v e_{v'} M_U(2).
\end{multline*}
Similarly,
\begin{multline*}
n_K^2 \left(\mathbb{E}\sum_{v\in G_k} \xi_{vi}\xi_{vj} \right)^2 =\left(\sum_{v\in G_k} x_{vi} x_{vj} \mathbb{E}(Y_v-e_v)^2\right)^2 =\left(\sum_{v\in G_k} x_{vi} x_{vj} e_v\right)^2\\
=\sum_{v\in G_k} x_{vi}^2 x_{vj}^2 e_v^2 +2\sum_{v<v'\in G_k} x_{vi}x_{vj}x_{v'i}x_{v'j} e_v e_{v'}.
\end{multline*}
Combining these last two gives the second variance term in $R_2$:
\begin{multline*}
n_K^2 \var \sum_{v\in G_k} \xi_{vi}\xi_{vj} = \sum_{v\in G_k} x_{vi}^2 x_{vj}^2 e_v [e_v(3 M_U(2)-1)+1]\\
+2\sum_{v<v'\in G_k} x_{vi}x_{vj}x_{v'i}x_{v'j} e_v e_{v'} (M_U(2)-1).
\end{multline*}
Applying this with $i=j$ gives the first variance term in $R_2$:
\begin{multline*}
n_K^2 \var \sum_{v\in G_k} \xi_{vj}^2  = \sum_{v\in G_k} x_{vj}^4 e_v [e_v(3 M_U(2)-1)+1]\\
+2\sum_{v<v'\in G_k} x_{vj}^2 x_{v'j}^2 e_v e_{v'} (M_U(2)-1).
\end{multline*}

Next we consider $R_3$ for which we write $\xi_i'=(Y_i'-e_i)x_i$, $i\in G_k$, with $Y_i'$ having the same distribution as $Y_i$ but independent of all else. Then 
\begin{multline}\label{pois.R3.crude}
\mathbb{E}\left(\sum_{j = 1}^d \abs{\xi'_{ij} - \xi_{ij}}\right)^3 = \mathbb{E}\left(n_K^{-1/2} \sum_{j = 1}^d |x_{ij}| |Y_i-Y_i'| \right)^3 \\
= n_K^{-3/2}\left( \sum_{j = 1}^d |x_{ij}|\right)^3 \mathbb{E}[|Y_i-Y_i'|^3],
\end{multline} 
and we crudely bound the last term by
\begin{multline}\label{pois.R3.crude2}
\mathbb{E}[|Y_i-Y_i'|^3]\le \mathbb{E}[(Y_i\vee Y_i')^3]\le \mathbb{E}[Y_i^3 +Y_i^{\prime 3}] = 2\mathbb{E}[Y_i^3]\\
=2[e_i + 3e_i^2 M_U(2)+ e_i^3 M_U(3)],
\end{multline}
using \eqref{pois.Y.mom} with ${3\brace 2}=3$ for the last equality.  Plugging this last into   \eqref{pois.R3.crude} gives 
\begin{equation*}
R_3\le 2 n_K^{-3/2} \sum_{i=1}^{n_K} e_i(1 + 3e_i M_U(2)+ e_i^2 M_U(3)) \left( \sum_{j = 1}^d |x_{ij}|\right)^3.
\end{equation*}

Lastly we address $R_1$. First we note that \eqref{R1.def.L2} vanishes: $\nabla S$ is constant by \eqref{pois.S} and thus has zero standard deviation, and the last term is zero by \eqref{pois.ES=0}. The  matrix~$\wtilde{I}_k$ was calculated in \eqref{pois.Itild.t0}, thus we find \eqref{R1.def.L3}. To find the functions~$M$ satisfying \eqref{M.fcn.bd}, taking an additional derivative of \eqref{pois.S} gives
\begin{equation*}
\abs{ \frac{\partial^3 \log f^{n_k}(y^{n_k}|\theta)}{\partial \theta_i \partial \theta_u \partial \theta_j}  }  \le \sum_{n=1}^{n_k} e^{x_n^T \theta}\left|x_{ni}x_{nu}x_{nj}\right|.
\end{equation*}
For $\eps>0$, if $|\theta_j - \theta_{0 j}| < \eps$ for all $j \in [d]$ then
\begin{equation*}
e^{x_n^T \theta} = e^{x_n^T \theta_0}e^{x_n^T (\theta-\theta_0)}\le e_n \exp\left(\eps \sum_{l=1}^d|x_{nl}|\right).
\end{equation*}
Thus taking
\begin{equation}\label{pois.M.fcn.def}
M^k_{iuj}(y^{n_k}) = \sum_{n=1}^{n_k} e_n \left|x_{ni}x_{nu}x_{nj}\right| \exp\left(\eps \sum_{l=1}^d|x_{nl}|\right)
\end{equation} satisfies \eqref{M.fcn.bd} and, being a constant function of $y^{n_k}$, the expectation in \eqref{EM.fcn.bd} is simply the square of \eqref{pois.M.fcn.def}.

We summarize our calculations for the Poisson GLMM in the following corollary.

\begin{corollary}\label{cor:pois}
Let $Y_1, Y_2, \dots, Y_{n_K}$ be observations from the Poisson GLMM \eqref{eq:poisson_glmm} satisfying \eqref{EU1}, the covariate vectors~$x_1,\ldots, x_{n_K}\in\mathbb{R}^d$ are full rank~$d$, and $\hat{\theta}^K$ maximizes  the misspecified likelihood~$f(\cdot|\theta)$ given by \eqref{pois.f.def}. For $Z \sim \mN(0, \mI_q)$, any $h \in \mathcal{H}$, and $\eps>0$,
\begin{multline}\label{pois.cor.bd}
\abs{\mathbb{E} [h(\sqrt{n_K} J^{-1/2} (\hat{\theta}^K - \theta_0^K))] - \mathbb{E} [ h(Z)]} \\
\le \frac{|h|_1}{\sqrt{n_K}} R_1 + \frac{ q^2 \tau^2 |h|_2}{4} R_2 + \frac{  q^3 \tau^3 |h|_3}{12} R_3 + 2|h|_0\frac{\mathbb{E}Q}{\eps},
\end{multline}
where $J$ is given by \eqref{eq:Jn_misspecified_definition} with $\theta_{n_K}^*=\theta_0$ and $H_k, I_k$ given by  \eqref{pois.H.t0} and \eqref{pois.I.t0}, $\theta_0^K=[\theta_0;\ldots; \theta_0]$, 
\begin{align*}
R_1 &=  \frac{1}{2}\sum_{k_1=1}^K \sum_{k_2=1\vee (k_1-1)}^{k_1} \sum_{i,j=1}^d  \left| \wtilde{I}_{k_1}(\theta_0)_{ij}^{-1/2}\right|  \\
&\times \sum_{l,l' = 1}^d \sqrt{\mathbb{E}[(Q_{k_2l}Q_{k_2l'})^2\wedge \eps^4]}  \sum_{n=1}^{n_{k_2}} e_n \left|x_{nj}x_{nl}x_{nl'}\right| \exp\left(\eps \sum_{w=1}^d|x_{nw}|\right),\\
R_2 &= \frac{1}{n_K}\sum_{k = 1}^{K} \left\{ \sum_{j = 1}^d \left[ \sum_{v\in G_k} x_{vj}^4 e_v [e_v(3 M_U(2)-1)+1]\right. \right.\\
&+\left. 2\sum_{v<v'\in G_k} x_{vj}^2 x_{v'j}^2 e_v e_{v'} \var(e^U)\right]^{1/2}\\
&+ \sum_{i<j} \left[\sum_{v\in G_k} x_{vi}^2 x_{vj}^2 e_v [e_v(3 M_U(2)-1)+1]\right.\\
&+\left.\left. 2\sum_{v<v'\in G_k} x_{vi}x_{vj}x_{v'i}x_{v'j} e_v e_{v'} \var(e^U)\right]^{1/2} \right\},\\
R_3 &= 2 n_K^{-3/2} \sum_{i=1}^{n_K} e_i(1 + 3e_i M_U(2)+ e_i^2 M_U(3)) \left( \sum_{j = 1}^d |x_{ij}|\right)^3,
\end{align*}
$\tau$, $Q_{k_i}$, and $Q$ are given by \eqref{non.asymp.c.def} and \eqref{Q.def} with $\theta_{n_K}^*=\theta_0$, the $e_i$ are given by \eqref{pois.ei.def}, and $M_U$ is the MGF of the random effects~$U_k$ in \eqref{eq:poisson_glmm}.
\end{corollary}

All quantities in the bound have been determined explicitly except for the expectations $\mathbb{E}[(Q_{k_2l}Q_{k_2l'})^2\wedge \eps^4]$ and $\mathbb{E}Q$, which can be estimated by the bootstrap. Under the asymptotic regime in
Section~\ref{sec:results_for_misspecified}, the bound in the corollary is of order $O(1/\sqrt{n_K})$ which follows from arguments similar to those for the correctly specified 1-parameter exponential family in \citet[][Section~4.3]{Aronowitz25} and orders of central moments of MLEs \citep{jose_jose_2001,peers_iqbal1985}. 

\section{Data and simulation example}\label{sec:data}

As an illustration, in this section we exhibit a data set arising from a clinical trial of epileptics carried out by \citet{Leppik85} that fits
the set up of this paper, and compute the bound for the Poisson GLMM in Corollary~\ref{cor:pois} in the context of that data set. We continue to use the notation for the Poisson GLMM from Section~\ref{sec:PoissGLM}. In that trial, patients suffering from seizures were randomized to receive either the antiepileptic drug \textit{progabide} or a placebo as an adjuvant to standard chemotherapy, and at each of $K=4$ successive post-randomization clinical visits, the number of epileptic seizures occurring over the previous 2 weeks was reported, yielding cumulative sample sizes
\begin{equation}\label{n.epil}
\bm{n}_{epil}:=(n_1,n_2,n_3,n_4) = (59, 118, 177, 236).
\end{equation}
The seizure count data~$Y_i$, along with covariates age (logarithm of age in years), treatment (indicator of the \textit{progabide} group), and baseline seizure rate (computed as the logarithm of one quarter of the $8$-week prerandomization seizure count) are available from the R package \texttt{MASS::epil} \citep{Venables02} or \citet[][Table~2]{Thall90}. With $K=4$ groups and $d=4$ covariates, $\what{\theta}^K$ is $q=16$ dimensional.

\begin{figure}
\centering
\includegraphics[width=.9\linewidth]{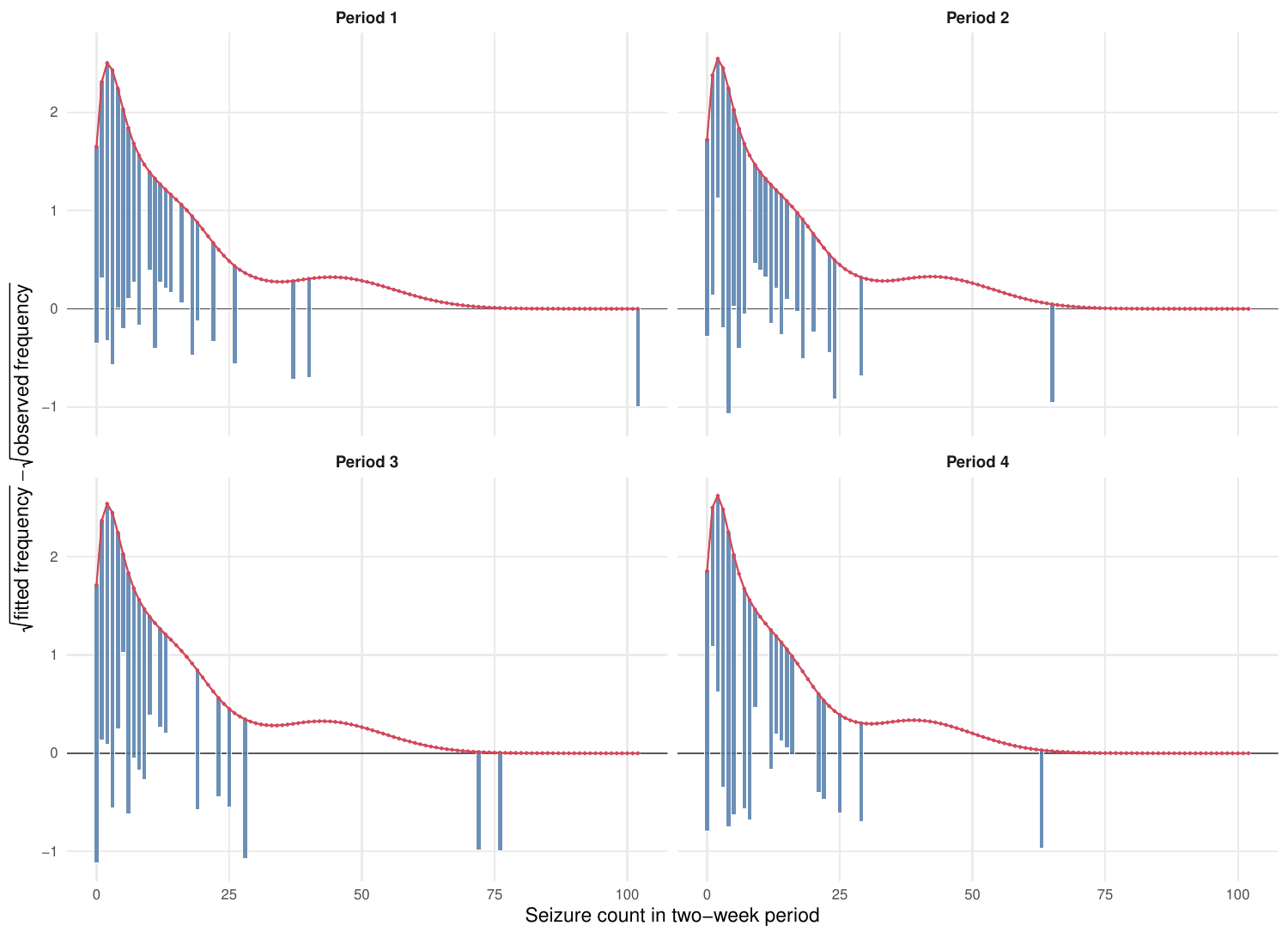}
\caption{Hanging Rootograms for epilepsy seizure counts by period: Fitted frequencies from Poisson GLMM of with covariates treatment, baseline seizure count, age, and random period effect. The full count range is shown.}
\label{fig:root}
\end{figure}

\begin{figure}
\centering
\includegraphics[width=0.9\linewidth]{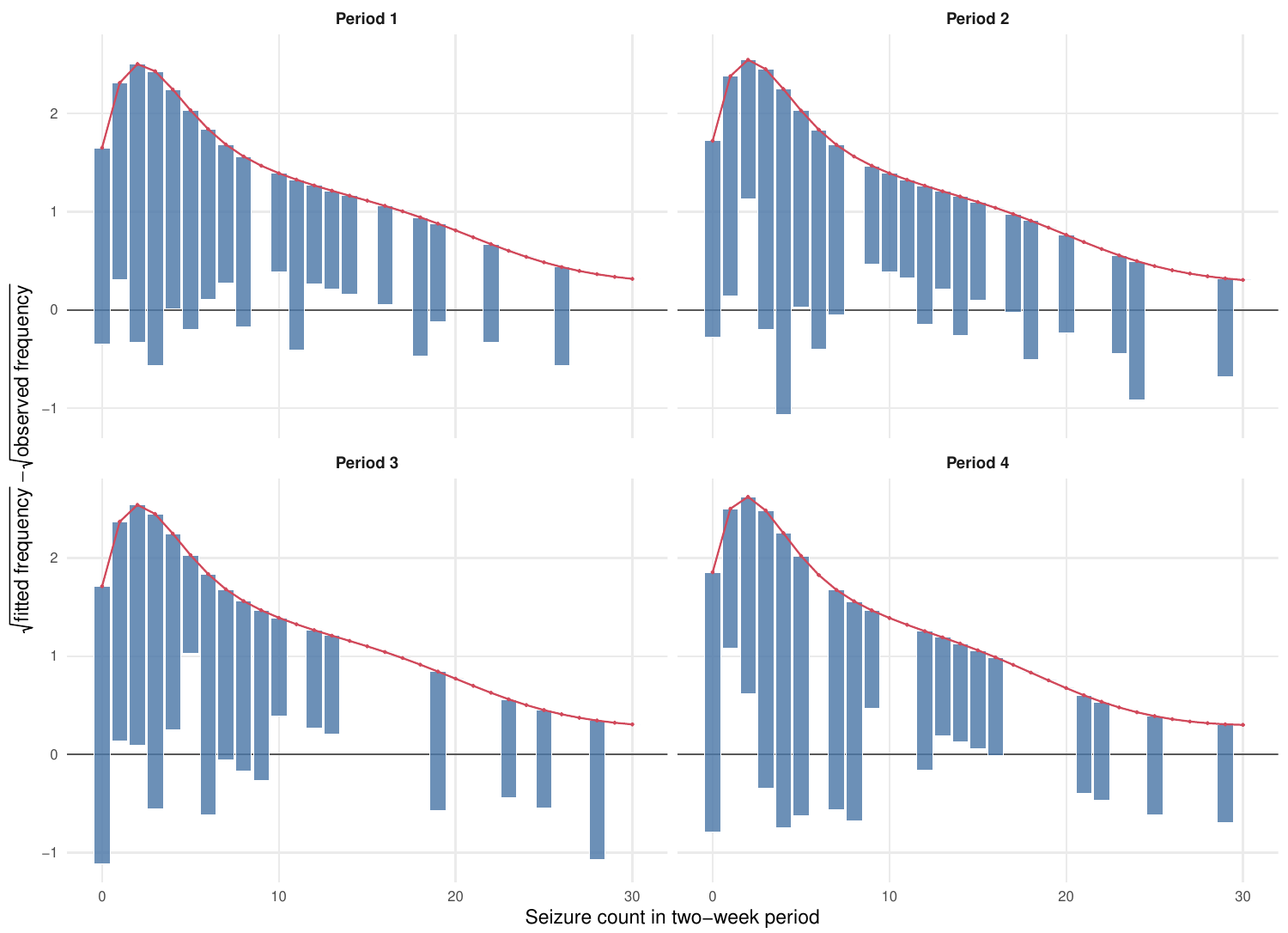}
\caption{Hanging Rootograms for epilepsy seizure counts by period: Fitted frequencies from Poisson GLMM of with covariates treatment, baseline seizure count, age, and random period effect. The figures are zoomed in to show count range $0$-$30$; the full range is shown in Figure~\ref{fig:root}.}
\label{fig:root.trunc}
\end{figure}

To visualize this data, Figure~\ref{fig:root} shows hanging rootograms \citep{Kleiber16} of the counts resulting from fitting a Poisson GLMM~\eqref{eq:poisson_glmm} to each period separately.  In Figure~\ref{fig:root} we see that the data includes a small number of counts exceeding $50$, so  Figure~\ref{fig:root.trunc} shows the same rootograms with the $x$-axes zoomed in to the $0$-$30$ range. In both figures, for visualization purposes, the data has been segregated by period, but the analysis which follows analyzes it cumulatively as a group sequential trial. In these figures and below, the period random effects are modeled as 
\begin{equation}\label{ex.U.dist}
U_k\sim\mN(-\sigma^2/2, \sigma^2),
\end{equation} for which $M_U(a)=\exp[ a(a-1)\sigma^2 /2]$ and thus satisfying \eqref{EU1}, where the parameter~$\sigma$ is estimated from the data. 

The seizure data summarized in Figures~\ref{fig:root} and \ref{fig:root.trunc} represents a single realization of the Poisson counts $Y_1,\ldots, Y_{n_K}$ and vector of estimates $\what{\theta}^K$ described in Section~\ref{sec:PoissGLM}. In order to investigate the behavior of the bound in Corollary~\ref{cor:pois} in this setting, we simulated data in the setting of this seizure count data in order to estimate the distributional distance on the left hand side of \eqref{pois.cor.bd}, as well as the corresponding upper bound on the right hand side. This was done by setting the true value of the regression coefficient vector~$\theta_0$ to be that calculated from the full seizure count data set, using the same number of $K=4$ groups, and then simulating random effects according to \eqref{ex.U.dist} and then Poisson counts according to \eqref{eq:poisson_glmm} with sample sizes 
\begin{equation}\label{sim.SS}
m\cdot \bm{n}_{epil}
\end{equation}
proportional to the actual sample sizes~\eqref{n.epil} for values of a multiplier $m=1,10,10^2,\ldots, 10^6$. For multiplier~$m$, each covariate vector was replicated $m$ times and the conditional Poisson distribution was used to generate responses.

Table~\ref{tab:epil.bd} contains the bound~\eqref{pois.cor.bd} in Corollary~\ref{cor:pois} broken down by the four terms there, calculated according to those quantities' definitions, and their percentage of the total bound.  The quantities involving expectations of $Q_{kl}$ in $R_1$, and of $Q$ in the last term of \eqref{pois.cor.bd}, were estimated using the parametric bootstrap with $500$ replications per table row.  Finally, since the first and last term on the right hand side of \eqref{pois.cor.bd} involve the free parameter~$\eps>0$, the sum of these two terms was numerically optimized by varying $\eps$, and the value of $\eps$ used for each row is included in the last column of the table. All this was carried out using the Gaussian test function $h(z)= (2\pi)^{-q/2} \exp(-z^T z/2)$ for $z\in\mathbb{R}^q$, with $q=16$,  a common choice in studies of distributional distances \citep[e.g.][]{Cucker07, Smola04}, for which basic calculations give 
$|h|_0 = (2\pi)^{-q/2}$, $|h|_1 = (2\pi)^{-q/2} e^{-1/2}$, $|h|_2 = (2\pi)^{-q/2}$, and $|h|_3=(2\pi)^{-q/2} \sqrt{6w^*}e^{-w^*/2}\approx 1.38 (2\pi)^{-q/2}$ where $w^*=3-\sqrt{6}$. 

From the table we see that in this example the $R_3$ term dominates the bound.  This is due in part to the crude bound used in \eqref{pois.R3.crude2}, and could be mitigated by utilizing the left hand side of \eqref{pois.R3.crude2}, computed by Monte Carlo, rather than the right hand side.  However, for this example, we left it in terms of the easily-computed closed form. In each successive row of the table, $m$ and hence the sample size increases by a factor
of $10$, and the bound consistently decreases by a factor close to $1/\sqrt{10}\approx .32$, reflecting
the fact that the bound is asymptotically $1/\sqrt{n_K}$, as discussed following Corollary~\ref{cor:pois}. Although the values of
the multiplier~$m$ required to make the bound small are sizable, they are comparable to 
those found for the state-of-the-art smooth function distance bounds for correctly-specified MLEs in both fixed-sample \citep{mvn,Anastasiou20}
and group sequential \citep{Aronowitz25}
settings, which do not have the additional variation due to misspecified model present here.  Still, the resulting bound is quite conservative with the left hand side of \eqref{pois.cor.bd} being of order $10^{-10}$ in this case for the sample sizes considered; this conservativeness is common for explicit Stein bounds of this type, and the simulation is intended to illustrate the computability, component structure, and rate of the bound.

\begin{table}[ht]
\centering
\scriptsize
\begin{tabular}{cllllcc}
\toprule
$m$ & $R_1$ term & $R_2$ term & $R_3$ term & $\mathbb{E}Q$ term & Bound & $\eps$ \\
\midrule
$1$      & $2.43\times 10^{-6}$ (.0\%) & .00160 (.6\%) & .28881 (99.4\%) & $4.32\times 10^{-6}$ (.0\%) & .290413 & .02314 \\
$10$     & $2.51\times 10^{-6}$ (.0\%) & .00052 (.6\%) & .08752 (99.4\%) & $3.47\times 10^{-6}$ (.0\%) & .088047 & .01557 \\
$10^2$    & $3.13\times 10^{-6}$ (.0\%) & .00022 (.8\%) & .02710 (99.2\%) & $3.24\times 10^{-6}$ (.0\%) & .027335 & .01396 \\
$10^3$  & $4.04\times 10^{-6}$ (.0\%) & .00016 (1.9\%) & .00851 (98.0\%) & $4.05\times 10^{-6}$ (.0\%) & .008682 & .01148 \\
$10^4$  & $5.17\times 10^{-6}$ (.2\%) & .00015 (5.4\%) & .00269 (94.2\%) & $5.39\times 10^{-6}$ (.2\%) & .002850 & .00855 \\
$10^5$ & $6.60\times 10^{-6}$ (.7\%) & .00015 (15.5\%) & .00080 (83.1\%) & $7.06\times 10^{-6}$ (.7\%) & .000963 & .00631 \\
$10^6$ & $8.42\times 10^{-6}$ (2.2\%) & .00015 (39.2\%) & .00022 (58.2\%) & $9.26\times 10^{-6}$ (2.4\%) & .000377 & .00412 \\
\bottomrule
\end{tabular}
\caption{Components of the Corollary~\ref{cor:pois} bound~\eqref{pois.cor.bd} for the epilepsy data Poisson GLMM example and values of the sample size multiplier~$m$ in \eqref{sim.SS}. Percentages indicate the relative contribution of each term to the total bound, and $\eps$ is the optimizing value in the $R_1$ and $\mathbb{E}Q$ terms.}
\label{tab:epil.bd}
\end{table}

\section{Conclusions}\label{sec:conc}
The asymptotic normality of sequences of random variables has been a major topic of research since at least \citet{deMoivre38}. Despite this topic's rich history, non-asymptotic bounds to these limits have  only gained maturity in the last half-century or so, and use of such bounds as finite-sample approximations in statistics is even less fully developed. For MLEs the first finite-sample bound in a general setting was not given until 2017 by \cite{AR_bound}. The current paper extends this line of research to include group sequential analysis and dependence through misspecified models and quasi-MLE.


\begin{thebibliography}{}

\bibitem[Anastasiou, 2017]{mle_bound_m_dependence}
Anastasiou, A. (2017).
\newblock Bounds for the normal approximation of the maximum likelihood
  estimator from m-dependent random variables.
\newblock {\em Statistics \& Probability Letters}, 129:171--181.

\bibitem[Anastasiou, 2018]{mvn}
Anastasiou, A. (2018).
\newblock Assessing the multivariate normal approximation of the maximum
  likelihood estimator from high-dimensional, heterogeneous data.
\newblock {\em Electronic Journal of Statistics}, 12(2):3794--3828.

\bibitem[Anastasiou and Gaunt, 2020]{Anastasiou20}
Anastasiou, A. and Gaunt, R.~E. (2020).
\newblock Multivariate normal approximation of the maximum likelihood estimator
  via the delta method.
\newblock {\em Brazilian Journal of Probability and Statistics}, 34:136--149.

\bibitem[Anastasiou and Ley, 2017]{delta}
Anastasiou, A. and Ley, C. (2017).
\newblock Bounds for the asymptotic normality of the maximum likelihood
  estimator using the delta method.
\newblock {\em Latin American Journal of Probability and Statistics},
  14:153–171.

\bibitem[Anastasiou and Reinert, 2017]{AR_bound}
Anastasiou, A. and Reinert, G. (2017).
\newblock Bounds for the normal approximation of the maximum likelihood
  estimator.
\newblock {\em Bernoulli}, 23(1):191–218.

\bibitem[Aronowitz and Bartroff, 2025]{Aronowitz25}
Aronowitz, J. and Bartroff, J. (2025).
\newblock Finite-sample bounds to the normal limit under group sequential
  sampling.
\newblock {\em Brazilian Journal of Probability and Statistics}, 39(1):1 -- 18.

\bibitem[Bartroff et~al., 2013]{Bartroff13}
Bartroff, J., Lai, T.~L., and Shih, M. (2013).
\newblock {\em Sequential Experimentation in Clinical Trials: Design and
  Analysis}.
\newblock Springer, New York.

\bibitem[Berckmoes, 2018]{Berckmoes18b}
Berckmoes, B. (2018).
\newblock An approach theoretic version of {A}nscombe’s theorem with an
  application in biostatistics.
\newblock {\em Applied Categorical Structures}, 26(5):883--889.

\bibitem[Berckmoes et~al., 2018]{Berckmoes18}
Berckmoes, B., Ivanova, A., and Molenberghs, G. (2018).
\newblock On the sample mean after a group sequential trial.
\newblock {\em Computational Statistics \& Data Analysis}, 125:104--118.

\bibitem[Berckmoes et~al., 2020]{Berckmoes20}
Berckmoes, B., Ivanova, A., and Molenberghs, G. (2020).
\newblock On asymptotic normality in estimation after a group sequential trial.
\newblock {\em Sequential Analysis}, 39(4):443--466.

\bibitem[Breslow and Clayton, 1993]{Breslow93}
Breslow, N.~E. and Clayton, D.~G. (1993).
\newblock Approximate inference in generalized linear mixed models.
\newblock {\em Journal of the American Statistical Association}, 88(421):9--25.

\bibitem[Chen et~al., 2010]{larrybook}
Chen, L.~H., Goldstein, L., and Shao, Q.-M. (2010).
\newblock {\em Normal Approximation by {S}tein's Method}.
\newblock Springer, Berlin Heidelberg.

\bibitem[Cox and Hinkley, 2000]{Cox00}
Cox, D.~R. and Hinkley, D.~V. (2000).
\newblock {\em Theoretical Statistics}.
\newblock Chapman \& Hall Ltd.

\bibitem[Cucker and Zhou, 2007]{Cucker07}
Cucker, F. and Zhou, D.~X. (2007).
\newblock {\em Learning theory: an approximation theory viewpoint}, volume~24.
\newblock Cambridge University Press.

\bibitem[de~A.~Cysneiros et~al., 2001]{jose_jose_2001}
de~A.~Cysneiros, F.~J., dos Santos, S. J.~P., and Cordeiro, G.~M. (2001).
\newblock Skewness and kurtosis for maximum likelihood estimator in
  one-parameter exponential family models.
\newblock {\em Brazilian Journal of Probability and Statistics}, 15(1):85--105.

\bibitem[de~Moivre, 1738]{deMoivre38}
de~Moivre, A. (1738).
\newblock {\em The Doctrine of Chances: Or, a Method of Calculating the
  Probabilities of Events in Play}.
\newblock A. Millar, London, 2nd edition.

\bibitem[Gaunt, 2016]{Gaunt16}
Gaunt, R.~E. (2016).
\newblock Rates of convergence in normal approximation under moment conditions
  via new bounds on solutions of the {S}tein equation.
\newblock {\em Journal of Theoretical Probability}, 29(1):231--247.

\bibitem[Gaunt and Li, 2023]{Gaunt23}
Gaunt, R.~E. and Li, S. (2023).
\newblock Bounding {K}olmogorov distances through {W}asserstein and related
  integral probability metrics.
\newblock {\em Journal of Mathematical Analysis and Applications},
  522(1):126985.

\bibitem[Greene, 2008]{Greene08}
Greene, W. (2008).
\newblock {\em Econometric Analysis}.
\newblock Pearson/Prentice Hall, United Kingdom.

\bibitem[Hayashi, 2011]{Hayashi11}
Hayashi, F. (2011).
\newblock {\em Econometrics}.
\newblock Princeton University Press, Princeton, New Jersey.

\bibitem[Hoadley, 1971]{hoadley1971}
Hoadley, B. (1971).
\newblock Asymptotic properties of maximum likelihood estimators for the
  independent not identically distributed case.
\newblock {\em Ann. Math. Statist.}, 42(6):1977--1991.

\bibitem[Jennison and Turnbull, 1997]{jt}
Jennison, C. and Turnbull, B.~W. (1997).
\newblock Group sequential analysis incorporating covariate information.
\newblock {\em Journal of the American Statistical Association},
  92(440):1330--1341.

\bibitem[Jennison and Turnbull, 2000]{sequential_book}
Jennison, C. and Turnbull, B.~W. (2000).
\newblock {\em Group Sequential Methods with Applications to Clinical Trials}.
\newblock Chapman \& Hall, Boca Raton, Florida.

\bibitem[Jiang, 1999]{glmm}
Jiang, J. (1999).
\newblock Conditional inference about generalized linear mixed models.
\newblock {\em Ann. Statist.}, 27(6):1974--2007.

\bibitem[Kim and Tsiatis, 2020]{Kim20}
Kim, K. and Tsiatis, A. (2020).
\newblock Independent increments in group sequential tests: a review.
\newblock {\em {SORT}: Statistics and operations research transactions},
  44(2):0223--264.

\bibitem[Kleiber and Zeileis, 2016]{Kleiber16}
Kleiber, C. and Zeileis, A. (2016).
\newblock Visualizing count data regressions using rootograms.
\newblock {\em The American Statistician}, 70(3):296--303.

\bibitem[Kohavi et~al., 2020]{Kohavi20}
Kohavi, R., Tang, D., and Xu, Y. (2020).
\newblock {\em Trustworthy online controlled experiments: A practical guide to
  {A}/{B} testing}.
\newblock Cambridge University Press, Cambridge.

\bibitem[Leppik et~al., 1985]{Leppik85}
Leppik, I.~E., Dreifuss, F.~E., Bowman, T., Santilli, N., Jacobs, M., Crosby,
  C., Cloyd, J., Stockman, J., Graves, N., Sutula, T., et~al. (1985).
\newblock A double-blind crossover evaluation of progabide in partial seizures.
\newblock {\em Neurology}, 35(4):285.

\bibitem[Lin et~al., 2022]{Lin22}
Lin, S., Scholtens, D., and Datta, S. (2022).
\newblock {\em Bioinformatics Methods: From Omics to Next Generation
  Sequencing}.
\newblock Chapman and Hall/CRC.

\bibitem[McCullagh, 1983]{mccullagh1983}
McCullagh, P. (1983).
\newblock Quasi-likelihood functions.
\newblock {\em Ann. Statist.}, 11(1):59--67.

\bibitem[McCulloch et~al., 2011]{mcculloch2011_book}
McCulloch, C., Searle, S., and Neuhaus, J. (2011).
\newblock {\em Generalized, Linear, and Mixed Models}.
\newblock Wiley Series in Probability and Statistics. Wiley.

\bibitem[Mc{C}ulloch and Neuhaus, 2012]{McCulloch12}
Mc{C}ulloch, C.~E. and Neuhaus, J.~M. (2012).
\newblock Generalized linear mixed models.
\newblock In {\em Encyclopedia of Environmetrics}. John Wiley \& Sons, Ltd, 2nd
  edition.

\bibitem[Parast and Bartroff, 2024]{Parast24}
Parast, L. and Bartroff, J. (2024).
\newblock Group sequential testing of a treatment effect using a surrogate
  marker.
\newblock {\em Biometrics}, 80(4):ujae108.

\bibitem[Peers and Iqbal, 1985]{peers_iqbal1985}
Peers, H.~W. and Iqbal, M. (1985).
\newblock Asymptotic expansions for confidence limits in the presence of
  nuisance parameters, with applications.
\newblock {\em Journal of the Royal Statistical Society. Series B
  (Methodological)}, 47(3):547--554.

\bibitem[Pinelis, 2017]{kol_mle_bound}
Pinelis, I. (2017).
\newblock Optimal-order uniform and nonuniform bounds on the rate of
  convergence to normality for maximum likelihood estimators.
\newblock {\em Electronic Journal of Statistics}, 11:1160--1179.

\bibitem[Powell, 2022]{Powell22}
Powell, W.~B. (2022).
\newblock {\em Reinforcement Learning and Stochastic Optimization}.
\newblock John Wiley \& Sons, Ltd, Princeton, NJ.

\bibitem[Reinert and R\"{o}llin, 2009]{RR_exchangeable_pair_bound}
Reinert, G. and R\"{o}llin, A. (2009).
\newblock Multivariate normal approximation with {S}tein's method of
  exchangable pairs under a general linearity condition.
\newblock {\em The Annals of Probability}, 36(6):2150--2173.

\bibitem[Riordan, 1937]{Riordan37}
Riordan, J. (1937).
\newblock Moment recurrence relations for binomial, {P}oisson and
  hypergeometric frequency distributions.
\newblock {\em The Annals of Mathematical Statistics}, 8(2):103--111.

\bibitem[Smola and Sch{\"o}lkopf, 2004]{Smola04}
Smola, A.~J. and Sch{\"o}lkopf, B. (2004).
\newblock A tutorial on support vector regression.
\newblock {\em Statistics and Computing}, 14:199--222.

\bibitem[Spiessens et~al., 2000]{Spiessens00}
Spiessens, B., Lesaffre, E., Verbeke, G., Kim, K., and DeMets, D.~L. (2000).
\newblock An overview of group sequential methods in longitudinal clinical
  trials.
\newblock {\em Statistical Methods in Medical Research}, 9(5):497--515.

\bibitem[Stein, 1972]{Stein72}
Stein, C. (1972).
\newblock A bound for the error in the normal approximation to the distribution
  of a sum of dependent random variables.
\newblock In {\em Proc.\ Sixth Berkeley Symp.\ Math.\ Stat.\ Prob.}, pages
  583--602. Univ.\ of California Press.

\bibitem[Thall and Vail, 1990]{Thall90}
Thall, P.~F. and Vail, S.~C. (1990).
\newblock Some covariance models for longitudinal count data with
  overdispersion.
\newblock {\em Biometrics}, 46(3):657--671.

\bibitem[Ulyanov, 1979]{ulyanov1979}
Ulyanov, V. (1979).
\newblock On more precise convergence rate estimates in the central limit
  theorem.
\newblock {\em Theory of Probability \& Its Applications}, 23(3):660--663.

\bibitem[Ulyanov, 1986]{ulyanov1986}
Ulyanov, V. (1986).
\newblock Normal approximation for sums of non-identically distributed random
  variables in {H}ilbert spaces.
\newblock {\em Acta Scientiarum Mathematicarum}, 50:411--419.

\bibitem[Ulyanov, 1987]{ulyanov1987}
Ulyanov, V. (1987).
\newblock Asymptotic expansions for distributions of sums of independent random
  variables in {H}.
\newblock {\em Theory of Probability \& Its Applications}, 31(1):25--39.

\bibitem[van Lint and Wilson, 2001]{vanLint01}
van Lint, J.~H. and Wilson, R.~M. (2001).
\newblock {\em A Course in Combinatorics}.
\newblock Cambridge University Press, 2nd edition.

\bibitem[Venables and Ripley, 2002]{Venables02}
Venables, W.~N. and Ripley, B.~D. (2002).
\newblock {\em Modern Applied Statistics with {S}}.
\newblock Springer, New York, 4th edition.
\newblock ISBN 0-387-95457-0.

\bibitem[Wang et~al., 2022]{Wang22}
Wang, T., Graves, B., Rosseel, Y., and Merkle, E.~C. (2022).
\newblock Computation and application of generalized linear mixed model
  derivatives using lme4.
\newblock {\em Psychometrika}, 87(3):1173--1193.

\bibitem[White, 1982]{white1982}
White, H. (1982).
\newblock Maximum likelihood estimation of misspecified models.
\newblock {\em Econometrica}, 50(1):1--25.

\bibitem[White, 1996]{white1996_book}
White, H. (1996).
\newblock {\em Estimation, Inference and Specification Analysis}.
\newblock Econometric Society Monographs. Cambridge University Press.

\end{thebibliography}

\def\cprime{$'$}

\appendix
\section{Regularity conditions}\label{sec:reg.conds}

The regularity conditions~\ref{reg:1st}-\ref{reg:CLT}
below are assumed for White's \citeyearpar{white1996_book} fixed sample result~\eqref{white.CLT} and our group sequential generalization of it in Theorem~\ref{thm:qmle.seq.asymp.norm}.  These conditions are a similar to conditions needed for the limiting normality of MLEs in correctly specific models in the fixed-sample \citep[see][]{hoadley1971} and group sequential settings  \citep[see][]{Aronowitz25,jt}, the difference being conditions concerning convergence of the estimator to the true parameter value in the latter setting, which are replaced by other arguments in our misspecified model setting. We remind the reader the $t$ denotes the dimension of each observation~$Y_n$. Unless stated otherwise, the following conditions are assumed to hold for all $\theta\in\Theta$ and $n=1,2,\ldots$.

\begin{enumerate}[label=(C\arabic*)]
\item\label{reg:1st} $\Theta$ is a compact subset of $\mathbb{R}^d$, $f_n(\cdot|\cdot, \theta): \mathbb{R}^{tn}\To \mathbb{R}^+$ is measurable,  and $f_n(Y_n|Y^{n-1},\cdot)$ is a.s.\ continuously second order differentiable on $\Theta$. 
\item\label{reg:LLN} $\E{\nabla^j\log f_n(Y_n|Y^{n-1}, \theta)}$ exists and $\{\nabla^j\log f_n(Y_n|Y^{n-1},\theta)\}$ obeys a weak LLN, for $j=0,1,2$. 
\item Differentiation of $\log f^n(Y^n| \theta)$ up to second order with respect to $\theta$ may be interchanged with expectation, i.e.,
\begin{equation*}
\nabla^j \E{\log f^n(Y^n| \theta)} = \E{\nabla^j \log f^n(Y^n| \theta)}<\infty\qmq{for}j=1,2, 
\end{equation*} 
and these are continuous on $\Theta$ uniformly in $n$. 
\item $\{(1/n) \log f^n(Y^n| \theta)\}$ is  $O(1)$ uniformly on $\Theta$ and, for all $n$, the maximizer~$\theta_n^*$ of $\E{\log f^n(Y^n| \theta)}$ is unique and  in the interior of $\Theta$.  
\item $\{\nabla^2 \E{\log f^n(Y^n| \theta_n^*)}\}$ is $O(1)$ and negative definite, uniformly in $n$. 
\item\label{reg:CLT} $\{n^{-1/2} \nabla  \log f^n(Y^n| \theta_n^*)\}$ obeys a central limit theorem and has a covariance matrix that is positive definite, uniformly in $n$, and $O(1)$. 
\end{enumerate}

The LLN and CLT conditions in \ref{reg:LLN} and \ref{reg:CLT}, respectively, are assumed in that form for generality; \citet[][pp.~27, 91]{white1996_book} lists a variety of sufficient conditions on the $f_n$ and $Y^n$ for these.






\section{Proof of Theorem~\ref{thm:qmle.seq.asymp.norm}}\label{sec:GS.asymp.pf} 
With the notation \eqref{S.set.notation}-\eqref{Gk.def}, for brevity we will denote the $\mathbb{R}^q$ column vectors
\begin{align}
S&=[S([n_1], \theta_{n_K}^*); S([n_2], \theta_{n_K}^*); \dots; S([n_K], \theta_{n_K}^*)],\label{S.vec.def}\\
S_G&=[S(G_1, \theta_{n_K}^*); S(G_2, \theta_{n_K}^*); \dots; S(G_K, \theta_{n_K}^*)].\label{SG.vec.def}
\end{align}
By applying the regularity conditions to the score function~$S(G_k,\theta_{n_K}^*)$ of each group and using independence between groups, we have 
\begin{equation*}
n_K^{-1/2}\wtilde{I}^{-1/2}S_G \overset{\mathcal{D}}{\longrightarrow} \mN(0, \mI_q)
\end{equation*} 
where 
\begin{equation}\label{D1.def}
\wtilde{I}=\diag{\wtilde{I}_1(\theta_{n_K}^*),\ldots, \wtilde{I}_K(\theta_{n_K}^*)}, 
\end{equation} which we use to denote a block-diagonal matrix with $\wtilde{I}_k$ given by \eqref{Itilde.def}. Since $S=AS_G$ for the linear transformation~$A$ with block structure
\begin{equation}\label{A.mat.def}
A_{[j][k]}= \mI_d 1\{j\ge k\},     
\end{equation}
$S=AS_G$ is also asymptotically multivariate normal:
\begin{equation}\label{S.th.star.Z}
n_K^{-1/2}\wtilde{I}^{-1/2}A^{-1}S \overset{\mathcal{D}}{\longrightarrow}  \mN(0,  \mI_q).
\end{equation}
The Taylor series expansion of $(1/n_k)S([n_k],\what{\theta}_{n_k})$ about $\theta_{n_K}^*$ is
\begin{equation*}
n_k^{-1} S([n_k],\what{\theta}_{n_k}) = n_k^{-1}S([n_k],\theta_{n_K}^*) + n_k^{-1}\nabla S([n_k],\theta^{\dagger}_{n_k})(\what{\theta}_{n_k} - \theta_{n_K}^*),
\end{equation*}
with $\theta^{\dagger}_{n_k}$ on the line segment between $\what{\theta}_{n_k}$ and $\theta_{n_K}^*$. Since $S([n_k], \what{\theta}_{n_k}) = 0$, this last becomes 
\begin{equation}\label{eq:taylor_series_Sk_misspecified_simplified}
n_k^{-1}S([n_k], \theta_{n_K}^*) = -n_k^{-1}\nabla S([n_k], \theta^{\dagger}_{n_k}) (\what{\theta}_{n_k} - \theta_{n_K}^*).
\end{equation}
The regularity conditions imply that 
$$\frac{n_k}{n_K}H_k(\theta)^{-1}\left(-n_k^{-1}\nabla S([n_k], \theta)\right)\To \mI_d.$$
Since $\what{\theta}_{n_k}-\theta_{n_K}^*$ is weakly consistent for $0$, so is $\theta^{\dagger}_{n_k}-\theta_{n_K}^*$, giving 
\begin{equation}\label{eq:avg_score_of_theta_dagger_limit}
\frac{n_k}{n_K}H_k(\theta_{n_K}^*)^{-1}\left(-n_k^{-1}\nabla S([n_k], \theta^{\dagger}_{n_k})\right)\To \mI_d.
\end{equation}
Letting 
\begin{equation}\label{D2.def}
H=\diag{H_1(\theta_{n_K}^*),\ldots, H_K(\theta_{n_K}^*)} 
\end{equation} and using \eqref{S.th.star.Z}, \eqref{eq:taylor_series_Sk_misspecified_simplified}, and \eqref{eq:avg_score_of_theta_dagger_limit}, we have
\begin{equation*}
\sqrt{n_K}\wtilde{I}^{-1/2}A^{-1}H (\what{\theta}^K - \theta^{*K}) \overset{\mathcal{D}}{\longrightarrow}
\mN(0, \mI_q).
\end{equation*}
All that remains is to show that 
\begin{equation}\label{asymp.J.deriv}
J=(H^{-1}A\wtilde{I}^{1/2})(H^{-1}A\wtilde{I}^{1/2})^T=H^{-1}(A\wtilde{I}A^T) (H^{-1})^T
\end{equation}
is given by the matrix~\eqref{eq:Jn_misspecified_definition}; the second equality in \eqref{asymp.J.deriv} uses that $\wtilde{I}$ is symmetric, being block-diagonal of covariance matrices. The desired equality then follows easily from matrix multiplication, using that $A\wtilde{I}A^T$ has block structure
\begin{equation}\label{AD1A.eq.Jtild}
(A\wtilde{I}A^T)_{[j][k]} = \sum_{i=1}^{j\wedge k} \wtilde{I}_i(\theta_{n_K}^*) = I_{j\wedge k}(\theta_{n_K}^*),\end{equation}
again using independence between groups, 
and $$H^{-1}=\diag{H_1(\theta_{n_K}^*)^{-1}, \ldots, H_K(\theta_{n_K}^*)^{-1}}$$ since $H$ is block diagonal. This completes the proof.\qed

\section{Proof of Theorem~\ref{thm:qmle_non_asymptotic}}\label{sec:qmle_non_asymptotic}

We continue to use the notation \eqref{S.vec.def}-\eqref{SG.vec.def}, and set
\begin{equation*}
\mu=\mathbb{E}S_G = [\mathbb{E}S(G_1, \theta_{n_K}^*); \mathbb{E}S(G_2, \theta_{n_K}^*); \dots; \mathbb{E}S(G_K, \theta_{n_K}^*)].
\end{equation*}
With $\wtilde{J}$ being the $q\times q$ matrix 
\begin{equation}\label{Jtild.def}
\wtilde{J}_{[j][k]}=I_{j\wedge k}(\theta_{n_K}^*),\quad j,k\in[K],
\end{equation}
and $A$ given  by \eqref{A.mat.def}, we will use the triangle inequality to bound
\begin{gather}
\abs{\mathbb{E} h(\sqrt{n_K} J^{-1/2} (\what{\theta}^K - \theta^{*K})) - \mathbb{E} h(Z)}\nonumber\\
\le\abs {\mathbb{E} h(n_K^{-1/2} \wtilde{J}^{-1/2} (S-A\mu)) - \mathbb{E} h(Z)} \label{eq:qmle_main triangle inequality 1st term} \\
+\abs{\mathbb{E} h(\sqrt{n_K} J^{-1/2} (\what{\theta}^K - \theta^{*K})) - \mathbb{E} h(n_K^{-1/2} \wtilde{J}^{-1/2} (S-A\mu))}. \label{eq:qmle_main triangle inequality 2nd term}
\end{gather}

\subsection{An upper bound for \eqref{eq:qmle_main triangle inequality 1st term}}
Set $W = n_K^{-1/2} (S_G-\mu)$ so $\mathbb{E}W=0$ and $\Var{W}=\wtilde{I}$ given by \eqref{D1.def}. Define the function $\tilde{h}(w) = h(\wtilde{J}^{-1/2} Aw)$ so the two terms in \eqref{eq:qmle_main triangle inequality 1st term} can be written 
\begin{equation}\label{hZ.nonasym}
h( n_K^{-1/2} \wtilde{J}^{-1/2} (S-A\mu)) = \tilde{h}(W) \qmq{and} h(Z)=h(\wtilde{J}^{-1/2} A \wtilde{I}^{1/2}Z)=\tilde{h}(\wtilde{I}^{1/2}Z),
\end{equation}
where this first equality for $h(Z)$ in \eqref{hZ.nonasym} holds because
\begin{equation}\label{eq:qmle_Sigma_to_the_neg_1/2}
\wtilde{J}^{-1/2} A \wtilde{I}^{1/2} = \mI_q,   
\end{equation} which follows from \eqref{AD1A.eq.Jtild}. Thus the term in \eqref{eq:qmle_main triangle inequality 1st term} is equal to $|\mathbb{E} \tilde{h}(W) - \mathbb{E} \tilde{h}(\wtilde{I}^{1/2}Z)|$, 
to which we apply \citet[][Theorem~2.1]{RR_exchangeable_pair_bound}, recorded as Theorem~\ref{RR_exchangeable_pair_bound} in Appendix~\ref{app:Stein.bg} for the reader's reference.  

To produce an exchangeable pair~$W'$ for $W$, let $i^*$ be a random index chosen uniformly from $[n_K]$, independent of all else.  Recalling \eqref{xi.def}, define $W'$ block-wise as
\begin{equation*}
W_{[k]}' = \begin{cases}
W_{[k]} - \xi_{i^*}+\xi_{i^*}',&\mbox{if $i^*\in G_k$,}\\
W_{[k]},&\mbox{otherwise}
\end{cases}
\end{equation*} for $k\in[K]$. We have
\begin{align*}
\mathbb{E} [W_{[k]}' - W_{[k]} | W] &= \sum_{i\in [n_K]}\mathbb{E} [W_{[k]}' - W_{[k]} | W, i^*=i] \mathbb{P}(i^*=i) \\
&=  \sum_{i\in G_k} \mathbb{E} [\xi_{i^*}' - \xi_{i^*} | W, i^*=i] \cdot \frac{1}{n_K}\\
&= n_K^{-1} \sum_{i\in G_k} \mathbb{E} [\xi_i'-\xi_i | W] \\
&= n_K^{-1} \sum_{i\in G_k} (\mathbb{E}[\xi_i']-\mathbb{E}[\xi_i|W]) \\
&= n_K^{-1} \sum_{i\in G_k} (\mathbb{E}[n_K^{-1/2}S(i,\theta_{n_K}^*)]-\mathbb{E}[n_K^{-1/2}S(i,\theta_{n_K}^*)|W]) \\
&= n_K^{-3/2} \left(\mathbb{E} \sum_{i\in G_k}  S(i, \theta_{n_K}^*) - \mathbb{E}\left[\left. \sum_{i\in G_k} S(i, \theta_{n_K}^*) \right| W\right]\right) \\
&= n_K^{-3/2} \left( \mathbb{E} S(G_k, \theta_{n_K}^*) - \mathbb{E}[S(G_k, \theta_{n_K}^*) | W]\right) \\
&= n_K^{-3/2} \left( \mu_{[k]} - (n_K^{1/2}W_{[k]}+\mu_{[k]})\right) \\
	&= -n_K^{-1} W_{[k]}.
\end{align*}
Applying Theorem~\ref{RR_exchangeable_pair_bound} with $\Lambda = (1/n_K)\mI_q$ and $R = 0$ and writing the result block-wise gives
\begin{align}
&\abs{\mathbb{E} \tilde{h}(W) - \mathbb{E} \tilde{h}(\wtilde{I}^{1/2}Z)}  \leq \nonumber\\
& \frac{n_K |\tilde{h}|_2}{4} \sum_{k_1, k_2 = 1}^{K}\sum_{i,j = 1}^{d} \sqrt{\Var{\mathbb{E} [(W_{[k_1]i}' - W_{[k_1]i})(W_{[k_2]j}' - W_{[k_2]j})| W]}} \label{eq:qmle_RR1}\\
&+\frac{n_K |\tilde{h}|_3}{12} \sum_{k_1, k_2, k_3 = 1}^{K}\sum_{i,j,u = 1}^{d} \mathbb{E} \left|(W_{[k_1]i}' - W_{[k_1]i}) (W_{[k_2]j}' - W_{[k_2]j})(W_{[k_3]u}' - W_{[k_3]u})\right| \label{eq:qmle_RR2}.
\end{align}
By definition of $W'$, the quantities inside the expectation in \eqref{eq:qmle_RR1} and \eqref{eq:qmle_RR2} vanish unless $k_1=k_2$ and $k_1=k_2=k_3$, respectively. Thus \eqref{eq:qmle_RR1} is equal to
\begin{equation*}
\frac{n_K |\tilde{h}|_2}{4} \sum_{k = 1}^{K}\sum_{i,j = 1}^{d} \sqrt{\Var {\mathbb{E} [(W_{[k]i}' - W_{[k]i})(W_{[k]j}' - W_{[k]j})| W]}} 
\end{equation*}
and \eqref{eq:qmle_RR2} is equal to
\begin{equation*}
\frac{n_K |\tilde{h}|_3}{12} \sum_{k = 1}^{K}\sum_{i,j,u = 1}^{d} \mathbb{E} \abs{(W_{[k]i}' - W_{[k]i})(W_{[k]j}' - W_{[k]j})(W_{[k]u}' - W_{[k]u})}. 
\end{equation*}
Next we write the $|\tilde{h}|_i$ in terms of the $|h|_i$. Using the change of variables $x = \wtilde{J}^{-1/2} Aw$, we have
\begin{equation*}
\frac{\partial \tilde{h}}{\partial w_j} = \sum_{i = 1}^{q} \frac{\partial h}{\partial x_i} \frac{\partial x_i}{\partial w_j} = \sum_{i = 1}^{q} \frac{\partial h}{\partial x_i} [\wtilde{J}^{-1/2} A]_{ij}
\end{equation*}
and
\begin{multline*}
\left|\frac{\partial^2 \tilde{h}}{\partial w_s \partial w_j}\right| = \left|\frac{\partial}{\partial w_s } \sum_{i_1 = 1}^{q} \frac{\partial h}{\partial x_{i_1}} [\wtilde{J}^{-1/2} A]_{i_1j}\right| = \left|\sum_{i_1 = 1}^{q} [\wtilde{J}^{-1/2} A]_{i_1j} \sum_{i_2 = 1}^{q} \frac{\partial^2 h}{\partial x_{i_2} \partial x_{i_1}} \frac{\partial x_{i_2}}{\partial w_s}\right|\\
= \left|\sum_{i_1 = 1}^{q} [\wtilde{J}^{-1/2} A]_{i_1j} \sum_{i_2 = 1}^{q} \frac{\partial^2 h}{\partial x_{i_2} \partial x_{i_1}} [\wtilde{J}^{-1/2} A]_{i_2s}\right|\\
\le \sum_{i_1, i_2 = 1}^q \left| [\wtilde{J}^{-1/2} A]_{i_1j}  [\wtilde{J}^{-1/2} A]_{i_2s}\right| \left | \frac{\partial^2 h}{\partial x_{i_2} \partial x_{i_1}} \right |  \le q^2 \tau^2 |h|_2,
\end{multline*}
where this last inequality uses \eqref{non.asymp.c.def} and \eqref{eq:qmle_Sigma_to_the_neg_1/2}. By similar arguments, 
\begin{equation*}
\left|\frac{\partial^3 \tilde{h}}{\partial w_u \partial w_s \partial w_j}\right| \le q^3 \tau^3 |h|_3,
\end{equation*}
and these inequalities give $|\tilde{h}|_2 \leq q^2 \tau^2 |h|_2$ and $|\tilde{h}|_3 \leq q^3 \tau^3 |h|_3$. Thus \eqref{eq:qmle_RR1} is $\le$
\begin{equation}
\frac{n_K q^2 \tau^2 |h|_2}{4} \sum_{k = 1}^{K} \sum_{i,j = 1}^{d} \bigg \{ \\ \Var{\mathbb{E} [(W_{[k]i}' - W_{[k]i})(W_{[k]j}' - W_{[k]j})\mid W]} \bigg \}^{1/2} \label{eq:qmle_RR1.1}
\end{equation}
and \eqref{eq:qmle_RR2} is $\le$
\begin{equation}\label{eq:qmle_RR2.1}
\frac{n_K q^3 \tau^3 |h|_3}{12} \sum_{k = 1}^{K} \sum_{i,j,u = 1}^{d} \mathbb{E} |(W_{[k]i}' - W_{[k]i})
(W_{[k]j}' - W_{[k]j})(W_{[k]u}' - W_{[k]u})|. 
\end{equation}
Next we bound the variance of the conditional expectations in \eqref{eq:qmle_RR1.1}. Let 
\begin{equation*}
C_2 = \frac{q^2 \tau^2 |h|_2}{4},\quad  C_3 = \frac{q^3 \tau^3 |h|_3}{12}.
\end{equation*}
With $\mathcal{A} = \sigma(Y^{n_K})$, the $\sigma$ algebra generated by the full data, since $\sigma(W) \subseteq \mathcal{A}$ we have
$\Var{ \mathbb{E} [\cdot | W]} \leq \Var{ \mathbb{E} [\cdot | \mathcal{A}]}$. Thus \eqref{eq:qmle_RR1.1} is $\le$
\begin{multline}
n_K C_2 \sum_{k = 1}^{K} \Big\{ \sum_{j = 1}^d \sqrt{\Var{\mathbb{E} [(\xi_{i^*j}' - \xi_{i^*j})^2 \mathbbm{1}\{i^* \in G_k \}| \mathcal{A}]}}\\
+ 2 \sum_{i < j} \sqrt{\Var{\mathbb{E} [(\xi_{i^*i}' - \xi_{i^*i})(\xi_{i^*j}' - \xi_{i^*j}) \mathbbm{1}\{i^* \in G_k \}| \mathcal{A}]}}\Big\}. \label{eq:qmle_square terms and cross terms}
\end{multline}
Using that $\mathbb{E} \xi_{ij}'=\mathbb{E} \xi_{ij}$ and $\mathbb{E}[\xi_{ij}^p|\mA] = \xi_{ij}^p$, the first variance term in \eqref{eq:qmle_square terms and cross terms} is
\begin{multline}
\Var{\mathbb{E} [(\xi_{i^*j}' - \xi_{i^*j})^2 \mathbbm{1}\{i^* \in G_k \}| \mathcal{A}]} = \var[ \mathbb{E} [(\xi_{i^*j}')^2 \mathbbm{1} \{i^* \in G_k \}]\\
-  2 n_K^{-1} \sum_{i \in G_k} \mathbb{E} [\xi_{ij}'] \mathbb{E} [\xi_{ij}| \mathcal{A}] + n_K^{-1} \sum_{i \in G_k} \mathbb{E} [\xi_{ij}^2| \mathcal{A}]]\\
=n_K^{-2}\cdot \var\sum_{i\in G_k}\xi_{ij}(\xi_{ij}-2\mathbb{E} \xi_{ij}).\label{R2.pt1}
\end{multline}
The second variance term in \eqref{eq:qmle_square terms and cross terms} is
\begin{multline}
\Var{\mathbb{E} [(\xi_{i^*i}' - \xi_{i^*i})(\xi_{i^*j}' - \xi_{i^*j}) \mathbbm{1}\{i^* \in G_k \}| \mathcal{A}]} = \var[\mathbb{E}[\xi_{i^*i}'\xi_{i^*j}' \mathbbm{1} \{i^* \in G_k \}]\\
- n_K^{-1} \sum_{v \in G_k} \mathbb{E} [\xi_{vi}'] \mathbb{E} [\xi_{vj}| \mathcal{A}] - n_K^{-1} \sum_{v \in G_k} \mathbb{E} [\xi_{vj}'] \mathbb{E} [\xi_{vi}| \mathcal{A}]  + n_K^{-1} \sum_{v \in G_k} \mathbb{E} [\xi_{vi} \xi_{vj}| \mathcal{A}]]\\
=n_K^{-2}\cdot\var \sum_{v \in G_k}\left(\xi_{vi}\xi_{vj}-\xi_{vj}\mathbb{E}\xi_{vi} - \xi_{vi}\mathbb{E}\xi_{vj}\right).\label{R2.pt2}
\end{multline}
Plugging \eqref{R2.pt1} and \eqref{R2.pt2} into \eqref{eq:qmle_square terms and cross terms}, we have that \eqref{eq:qmle_RR1.1} is $\le$
\begin{multline*}
C_2 \sum_{k = 1}^{K} \left\{ \sum_{j = 1}^d \left[\var\sum_{i\in G_k}\xi_{ij}(\xi_{ij}-2\mathbb{E} \xi_{ij})\right]^{1/2}\right. \\
+\left.  \sum_{i<j} \left[\var \sum_{v \in G_k}\left(\xi_{vi}\xi_{vj}-\xi_{vj}\mathbb{E}\xi_{vi} - \xi_{vi}\mathbb{E}\xi_{vj}\right)\right]^{1/2}\right\} = C_2R_2,
\end{multline*}
where $R_2$ is defined in the theorem. 

Returning to \eqref{eq:qmle_RR2.1}, it equals
\begin{align*}
&n_K C_3 \sum_{k = 1}^{K}\sum_{i,j,u = 1}^d  \mathbb{E} \left|(\xi'_{i^*i} - \xi_{i^*i})(\xi'_{i^*j} - \xi_{i^*j})(\xi'_{i^*u} - \xi_{i^*u}) \mathbbm{1}\{i^* \in G_k \}]\right|  \\
&=n_K C_3 \sum_{k = 1}^{K}\sum_{i,j,u = 1}^d  \frac{1}{n_K} \sum_{v \in G_k} \mathbb{E} \abs{(\xi'_{vi} - \xi_{vi})(\xi'_{vj} - \xi_{vj})(\xi'_{vu} - \xi_{vu})}   \\
&= C_3 \sum_{k = 1}^{K} \sum_{v \in G_k} \mathbb{E}  \sum_{i,j,u = 1}^d  \abs{\xi'_{vi} - \xi_{vi}} \abs{\xi'_{vj} - \xi_{vj}} \abs{\xi'_{vu} - \xi_{vu}}  \\
&= C_3 \sum_{k = 1}^{K} \sum_{v \in G_k} \mathbb{E} \left( \sum_{j = 1}^d \abs{\xi'_{vj} - \xi_{vj}}\right)^3  \\
&= C_3 \sum_{i = 1}^{n_K} \mathbb{E} \left( \sum_{j = 1}^d \abs{\xi'_{ij} - \xi_{ij}}\right)^3 \\
&= C_3 R_3.
\end{align*}
We have shown that \eqref{eq:qmle_main triangle inequality 1st term} is bounded above by $C_2 R_2+C_3 R_3$, the second and third terms in \eqref{eq:qmle_non_asymptotic_statement}.  Next we finish the proof by showing that \eqref{eq:qmle_main triangle inequality 2nd term} is bounded above by the remaining two terms in \eqref{eq:qmle_non_asymptotic_statement}.
	
\subsection{An upper bound for \eqref{eq:qmle_main triangle inequality 2nd term}} 

For a generic column vector argument
\begin{equation}\label{th.K.gen.def}
\theta^K=[\theta_1;\theta_2;\ldots; \theta_K]\in\mathbb{R}^q
\end{equation}
with $\theta_k\in\mathbb{R}^d$, further generalize \eqref{S.vec.def} by writing
\begin{equation}\label{S.th.K.def}
S(\theta^K) = [S([n_1], \theta_1); S([n_2], \theta_2); \dots; S([n_K], \theta_K)] \in\mathbb{R}^q.
\end{equation}
Below we will write $\theta_i^K\in\mathbb{R}$ for the $i$th entry of $\theta^K$, which should not be confused with the sub-vector $\theta_i\in\mathbb{R}^d$ in \eqref{th.K.gen.def}. Until this point we have used $\nabla$ to denote differentiation with respect to $\theta\in\mathbb{R}^d$, such as in \eqref{score.def} and \eqref{Hk.def}, but at this point we need to consider differentiation of $S(\theta^K)$ with respect to the full $q$-dimensional argument~$\theta^K$, which we denote by $\nabla_q$ for clarity. We note that, because of the structure~\eqref{S.th.K.def},  $\partial S(\theta^K)_j/\partial \theta_i^K=0$ for $i,j$ in different $d$-blocks of $[q]$.

A second order Taylor expansion of $S(\theta^K)$ about $\theta^{*K}$ evaluated at $\what{\theta}^K$ yields
\begin{multline}\label{S.Tay2.vec}
S(\what{\theta}^K) = S + \nabla_q S(\theta^{*K})(\what{\theta}^K - \theta^{*K})\\
+\frac{1}{2} \left[\begin{array}{c}
(\what{\theta}^K - \theta^{*K})^T \nabla_q^2 S(\theta_{\dagger}^K)_1 (\what{\theta}^K - \theta^{*K})\\
\vdots\\
(\what{\theta}^K - \theta^{*K})^T \nabla_q^2 S(\theta_{\dagger}^K)_q (\what{\theta}^K - \theta^{*K})
\end{array}\right],
\end{multline} where $\nabla_q^2 S(\theta_{\dagger}^K)_j$ denotes the $q\times q$ matrix with $(i,l)$ entry
$$\left[\nabla_q^2 S(\theta_{\dagger}^K)_j\right]_{il} = \left.\frac{\partial^2 S(\theta^K)_j}{\partial \theta_i^K \partial \theta_l^K}\right|_{\theta^K=\theta_{\dagger}^K}.$$
Let $V$ be the last vector in \eqref{S.Tay2.vec}.  With $H$ given by \eqref{D2.def}, adding $n_K H (\what{\theta}^K - \theta^{*K})$ to both sides of \eqref{S.Tay2.vec} and using that  $S(\what{\theta}^K)=0$ yields
\begin{equation}\label{eq:qmle_n_J_star}
n_K H (\what{\theta}^K - \theta^{*K}) = S +  (\nabla_q S(\theta^{*K})+n_K H)(\what{\theta}^K - \theta^{*K})+V.
\end{equation}
From \eqref{eq:Jn_misspecified_definition} and \eqref{Jtild.def} we have $\wtilde{J} = H J H$, hence $\wtilde{J}^{-1/2} = J^{-1/2} H^{-1}$. Thus multiplying both sides of \eqref{eq:qmle_n_J_star} by $n_K^{-1/2} \wtilde{J}^{-1/2}$ gives
\begin{multline} \label{eq:qmle_taylor form of (hat theta^K - theta^K)}
\sqrt{n_K} J^{-1/2} (\what{\theta}^K - \theta^{*K})  \\ 
= n_K^{-1/2} \wtilde{J}^{-1/2} \left[S +  (\nabla_q S(\theta^{*K})+n_K H)(\what{\theta}^K - \theta^{*K})+V\right].
\end{multline}
By \eqref{eq:qmle_Sigma_to_the_neg_1/2}, $\wtilde{J}^{-1/2}=\wtilde{I}^{-1/2}A^{-1}$ and, by \eqref{A.mat.def}, $A^{-1}$ has block structure
\begin{equation*}
(A^{-1})_{[j][k]}=\begin{cases}\mI_d,& j=k\\
-\mI_d,& j=k+1\\
0,&\mbox{otherwise,}     
\end{cases}
\end{equation*}
for $j,k\in[K]$, thus $\wtilde{J}^{-1/2}$ has block structure
\begin{equation}\label{J.tild.block}
(\wtilde{J}^{-1/2})_{[j][k]}=\begin{cases}\wtilde{I}_j(\theta_{n_K}^*)^{-1/2},& j=k\\
-\wtilde{I}_j(\theta_{n_K}^*)^{-1/2},& j=k+1\\
0,&\mbox{otherwise,}     
\end{cases}
\end{equation}
for $j,k\in[K]$. The following will be useful in relating $q$-dimensional objects to their $d$-block structure: If $B\in\mathbb{R}^{q\times q}$ and $v\in\mathbb{R}^q$ then
\begin{equation}\label{v.sum.block}
\sum_{i=1}^q (Bv)_i =  \sum_{k_1, k_2=1}^K \sum_{i,j=1}^d B_{[k_1][k_2]ij} v_{[k_2]j}.
\end{equation}

Letting
\begin{equation}\label{T3.def}
T_3 = h(n_K^{-1/2} \wtilde{J}^{-1/2} (S-A\mu+V)) - h(n_K^{-1/2} \wtilde{J}^{-1/2} (S-A\mu))
\end{equation}
and using a first order Taylor approximation and \eqref{eq:qmle_taylor form of (hat theta^K - theta^K)}, the quantity inside the expectation in \eqref{eq:qmle_main triangle inequality 2nd term} can be written 
\begin{multline*}
h \big (\sqrt{n_K} J^{-1/2} (\what{\theta}^K - \theta^{*K}) \big ) - h(n_K^{-1/2} \wtilde{J}^{-1/2} (S-A\mu))\\
=h \big (\sqrt{n_K} J^{-1/2} (\what{\theta}^K - \theta^{*K}) \big ) -h(n_K^{-1/2} \wtilde{J}^{-1/2} (S-A\mu+V)) +T_3\\
=\nabla_q h(x^{\dagger})\left( \sqrt{n_K} J^{-1/2} (\what{\theta}^K - \theta^{*K}) - n_K^{-1/2} \wtilde{J}^{-1/2} (S - A\mu+V)\right)+T_3\\
=\nabla_q h(x^{\dagger}) n_K^{-1/2} \wtilde{J}^{-1/2} \left((\nabla_q S(\theta^{*K})+n_K H)(\what{\theta}^K - \theta^{*K})+A\mu\right)+T_3\\
=:T_1+T_2+T_3.
\end{multline*}
Then \eqref{eq:qmle_main triangle inequality 2nd term} is equal to
\begin{equation*}
\left|E(T_1+T_2+T_3)\right|\le \left|ET_1\right|+ \left|ET_2\right|+ \left|ET_3\right|
\end{equation*} 
which we will bound from above by the sum of the first and last terms in \eqref{eq:qmle_non_asymptotic_statement}. For the remainder of the proof we denote the first two summations in \eqref{Rem1.def} by $\sum_{k_1-k_2\le 1}$.

For $T_1$, using \eqref{J.tild.block}, \eqref{v.sum.block}, and the block diagonal structure of $\nabla_q S(\theta^{*K})+n_K H$,  we have
\begin{multline}\label{T1.bd}
\frac{\sqrt{n_K} |\mathbb{E} T_1 |}{|h|_1} \leq \mathbb{E}\sum_{i=1}^q \left|\left[\wtilde{J}^{-1/2} (\nabla_q S(\theta^{*K})+n_K H)(\what{\theta}^K - \theta^{*K})\right]_i\right|\\
=\sum_{k_1, k_2=1}^K \sum_{i,j=1}^d  \left| \wtilde{J}_{[k_1][k_2]ij}^{-1/2}\right| \mathbb{E}\left|\left[(\nabla_q S(\theta^{*K})+n_K H)(\what{\theta}^K - \theta^{*K})\right]_{[k_2]j}\right|\\
=\sum_{k_1-k_2\le 1} \sum_{i,j=1}^d  \left| \wtilde{I}_{k_1}(\theta_{n_K}^*)_{ij}^{-1/2}\right| \mathbb{E}\left| \sum_{l=1}^d (\nabla_q S(\theta^{*K})+n_K H)_{[k_2][k_2]jl} (\what{\theta}^K - \theta^{*K})_{[k_2]l}\right|\\
\le \sum_{k_1-k_2\le 1} \sum_{i,j,l=1}^d  \left| \wtilde{I}_{k_1}(\theta_{n_K}^*)_{ij}^{-1/2}\right|  \mathbb{E}\left|(\nabla S([n_{k_2}],\theta_{n_K}^*)+n_K H_{k_2}(\theta_{n_K}^*))_{jl} (\what{\theta}_{n_{k_2}} - \theta_{n_K}^*)_l\right|\\
\le \sum_{k_1-k_2\le 1} \sum_{i,j,l=1}^d  \left| \wtilde{I}_{k_1}(\theta_{n_K}^*)_{ij}^{-1/2}\right|  \sqrt{\var[\nabla S([n_{k_2}],\theta_{n_K}^*)_{jl}]
\mathbb{E}[(\what{\theta}_{n_{k_2}} - \theta_{n_K}^*)_l]^2},
\end{multline}
where the final inequality uses  the Cauchy-Schwartz inequality and that $\mathbb{E}(\nabla S([n_{k_2}],\theta_{n_K}^*)) = -n_K H_{k_2}(\theta_{n_K}^*)$.

Similarly for $T_2$, using that $A\mu=\mathbb{E}S$, we have
\begin{equation}\label{T2.bd}
\frac{\sqrt{n_K} |\mathbb{E} T_2 |}{|h|_1} \le  \sum_{k_1-k_2\le 1} \sum_{i,j=1}^d  \left| \wtilde{I}_{k_1}(\theta_{n_K}^*)_{ij}^{-1/2}\right|  \left| \mathbb{E}S([n_{k_2}], \theta_{n_K}^*)_j \right|.\end{equation}

For $T_3$ we write
\begin{equation}\label{T3.1st.bd}
\mathbb{E}|T_3| = \mathbb{E}[\left.|T_3|\right| Q<\eps]\mathbb{P}(Q<\eps) + \mathbb{E}[\left.|T_3|\right| Q\ge \eps]\mathbb{P}(Q\ge \eps)
\end{equation} 
and we crudely bound the second term using \eqref{T3.def} and Markov's inequality: 
\begin{equation}\label{abs.T3.bd1}
\mathbb{E}[\left.|T_3|\right| Q\ge \eps]\mathbb{P}(Q\ge \eps) \le 2|h|_0\frac{\mathbb{E}Q}{\eps}.
\end{equation} 
For the first term in \eqref{T3.1st.bd}, we proceed similarly to $T_1$ and $T_2$ using the first order Taylor approximation
\begin{align*}
T_3 &= \nabla_q h(x^{\dagger\dagger})\left( \sqrt{n_K} J^{-1/2} (\what{\theta}^K - \theta^{*K}) - n_K^{-1/2} \wtilde{J}^{-1/2} (S - A\mu+V)\right) \\
&= \nabla_q h(x^{\dagger\dagger}) n_K^{-1/2} \wtilde{J}^{-1/2} V
\end{align*} so that
\begin{multline}\label{abs.T3.bd2}
\frac{2\sqrt{n_K}\; \mathbb{E}[\left.|T_3|\right| Q<\eps]}{|h|_1} \le
\sum_{k_1-k_2\le 1} \sum_{i,j=1}^d  \left| \wtilde{I}_{k_1}(\theta_{n_K}^*)_{ij}^{-1/2}\right| \mathbb{E}[\left.\left| 2V_{[k_2]j}\right|\right| Q<\eps]\\
\le \sum_{k_1-k_2\le 1} \sum_{i,j,l,l'=1}^d  \left| \wtilde{I}_{k_1}(\theta_{n_K}^*)_{ij}^{-1/2}\right| \\
\times \mathbb{E}\left[\left.\left|(\what{\theta}_{n_{k_2}}-\theta_{n_K}^*)_l (\what{\theta}_{n_{k_2}}-\theta_{n_K}^*)_{l'} \frac{\partial^2 S([n_{k_2}], \theta^\dagger)_j}{\partial \theta_l\theta_{l'}}\right|\right| Q<\eps\right],
\end{multline}
where here $\theta^\dagger\in \mathbb{R}^d$ lies between $\what{\theta}_{n_{k_2}}$ and $\theta_{n_K}^*$, and $\theta_l, \theta_{l'}$ denote entries of the generic  vector $\theta\in\mathbb{R}^d$ argument to $S([n_{k_2}],\cdot)$. This last term we bound using Cauchy-Schwarz and \eqref{M.fcn.bd}-\eqref{EM.fcn.bd}:
\begin{multline*}
\mathbb{E}\left[\left.\left|(\what{\theta}_{n_{k_2}}-\theta_{n_K}^*)_l (\what{\theta}_{n_{k_2}}-\theta_{n_K}^*)_{l'} \frac{\partial^2 S([n_{k_2}], \theta^\dagger)_j}{\partial \theta_l\theta_{l'}}\right|\right| Q<\eps\right] \\
\le \left\{\mathbb{E}[(Q_{k_2l}Q_{k_2l'})^2\wedge \eps^4] \; \mathbb{E} [ (  M_{jll'}^{k_2}(Y^{n_{k_2}}) )^2 | Q< \eps ]\right\}^{1/2}.
\end{multline*}
Plugging this last into \eqref{abs.T3.bd2}, and combining with \eqref{T1.bd}, \eqref{T2.bd}, \eqref{T3.1st.bd}, and \eqref{abs.T3.bd1} gives the desired result and completes the proof.\qed

\section{A result of \cite{RR_exchangeable_pair_bound}}\label{app:Stein.bg}

\begin{theorem}[\cite{RR_exchangeable_pair_bound}] \label{RR_exchangeable_pair_bound}
Assume that $(W, W')$ is an exchangeable pair of $\R^q$-valued random vectors such that 
\begin{equation*}
\mathbb{E} W = 0, \quad \mathbb{E} WW^T = \Sigma
\end{equation*}
with $\Sigma \in \R^{q \times q}$ symmetric and positive definite, and 
\begin{equation*} 
\mathbb{E} [W' - W \mid W] = -\Lambda W + R
\end{equation*}
for an invertible matrix $\Lambda$ and a $\sigma(W)$-measurable random vector $R$. Then, if $Z$ has $q$-dimensional standard normal distribution, we have for every three times differentiable function $h$, 
\begin{equation*}
\abs{\mathbb{E} h(W) - \mathbb{E} h (\Sigma^{1/2} Z)} \leq \frac{|h|_2}{4}A + \frac{|h|_3}{12}B + \Big (|h|_1 + \frac{1}{2} q \abs{\Sigma}^{1/2} |h|_2 \Big ) C,
\end{equation*}
where
\begin{align*}
A &= \sum_{i, j = 1}^{q} \lambda^{(i)} \sqrt{\Var{\mathbb{E} \{(W_i' - W_i)(W_j' - W_j) \mid W \}}}, \\
B &= \sum_{i, j, k = 1}^{q} \lambda^{(i)} \mathbb{E} \abs{(W_i' - W_i)(W_j' - W_j)(W_k' - W_k)}, \\
C &= \sum_{i,m = 1}^q \left| (\Lambda^{-1})_{mi}\right| \sqrt{\Var{R_i}}
\end{align*}
\end{theorem}

\end{document}